\documentclass[letterpaper]{amsart}
\pdfoutput=1
\newtheorem{theorem}{Theorem}[section]
\newtheorem{lemma}[theorem]{Lemma}

\newtheorem{proposition}[theorem]{Proposition}

\theoremstyle{definition}
\newtheorem{definition}[theorem]{Definition}

\theoremstyle{remark}
\newtheorem{remark}[theorem]{Remark}

\usepackage{graphicx}
\usepackage{color}
\usepackage{amsmath}
\usepackage{amsfonts}
\usepackage{amssymb}
\newcommand{\be}{\begin{equation}}

\newcommand{\ee}{\end{equation}}

\newcommand{\om}{\omega}

\newcommand{\clc}{{\stackrel{\scriptscriptstyle{LC}}{\nabla}}}
\newcommand{\lc}{{\stackrel{\scriptscriptstyle{LC}}{\Gamma}}\phantom{}}
\newcommand{\rilc}{{\stackrel{\scriptscriptstyle{LC}}{Ric}}\phantom{}}

\newcommand{\la}{\lambda}

\newcommand{\dz}{\wedge}

\newcommand{\R}{{\bf R }}

\newcommand{\ba}{\begin{array}}

\newcommand{\ea}{\end{array}}

\newcommand{\beq}{\begin{eqnarray}}

\newcommand{\eeq}{\end{eqnarray}}

\newtheorem{lm}{lemma}

\newtheorem{thee}{theorem}

\newtheorem{proo}{proposition}

\newtheorem{co}{corollary}

\newtheorem{rem}{remark}

\newtheorem{deff}{definition}

\newcommand{\bd}{\begin{deff}}

\newcommand{\ed}{\end{deff}}

\newcommand{\bl}{\begin{lm}}

\newcommand{\el}{\end{lm}}

\newcommand{\bp}{\begin{proo}}

\newcommand{\ep}{\end{proo}}

\newcommand{\bt}{\begin{thee}}

\newcommand{\et}{\end{thee}}

\newcommand{\bc}{\begin{co}}

\newcommand{\ec}{\end{co}}

\newcommand{\brm}{\begin{rem}}

\newcommand{\erm}{\end{rem}}

\newcommand{\der}{{\rm d}}

\usepackage{t1enc}
\def\frak{\mathfrak}

\newcommand{\newc}{\newcommand}

\renewcommand{\exp}{\operatorname{exp}}
\newcommand{\id}{\operatorname{id}}

\let\ccdot\cdot
\def\cdot{\hbox to 2.5pt{\hss$\ccdot$\hss}}

\newc{\aR}{\mbox{\boldmath{$ R$}}}
\newc{\aS}{\mbox{\boldmath{$ S$}}}
\newc{\aT}{\mbox{\boldmath{$ T$}}}
\newc{\aW}{\mbox{\boldmath{$ W$}}}

\newc{\aK}{\mbox{\boldmath{$ K$}}}
\newc{\aL}{\mbox{\boldmath{$ L$}}}

\newcommand{\bbC}{\mathbb{C}}

\usepackage{amssymb}
\usepackage{amscd}

\let\e=\varepsilon

\newcommand{\hook}{\raisebox{-0.35ex}{\makebox[0.6em][r]
{\scriptsize $-$}}\hspace{-0.15em}\raisebox{0.25ex}{\makebox[0.4em][l]{\tiny
 $|$}}}

\newcommand{\bas}{{\bf e}}

\newc{\obstrn}[2]{B^{#1}_{#2}}

\newcommand{\rpl}                         
{\mbox{$
\begin{picture}(12.7,8)(-.5,-1)
\put(0,0.2){$+$}
\put(4.2,2.8){\oval(8,8)[r]}
\end{picture}$}}

\newcommand{\lpl}                         
{\mbox{$
\begin{picture}(12.7,8)(-.5,-1)
\put(2,0.2){$+$}
\put(6.2,2.8){\oval(8,8)[l]}
\end{picture}$}}

\usepackage{ifthen}
\newcommand{\bbR}{\mathbb{R}}
\newcommand{\bbS}{\mathbb{S}}
\newcommand{\bbM}{\mathbb{M}}
\newcommand{\sping}{\mathbf{Spin}}
\newcommand{\spina}{\frak{spin}}
\newcommand{\sog}{\mathbf{SO}}

\newcommand{\slg}{\mathbf{SL}}
\newcommand{\glg}{\mathbf{GL}}

\newcommand{\soa}{\frak{so}}
\newcommand{\sla}{\frak{sl}}
\newcommand{\gla}{\frak{gl}}
\newcommand{\sua}{\frak{su}}

\newc{\tensor}[1]{#1}
\newc{\Mvariable}[1]{\mbox{#1}}
\newc{\down}[1]{{}_{#1}}
\newc{\up}[1]{{}^{#1}}

\newc{\JulyStrut}{\rule{0mm}{6mm}}
\newc{\midtenPan}{\mbox{\sf S}}
\newc{\midten}{\mbox{\sf T}}
\newc{\midtenEi}{\mbox{\sf U}}
\newc{\ATen}{\mbox{\sf E}}
\newc{\BTen}{\mbox{\sf F}}
\newc{\CTen}{\mbox{\sf G}}

\def\sideremark#1{\ifvmode\leavevmode\fi\vadjust{\vbox to0pt{\vss
 \hbox to 0pt{\hskip\hsize\hskip1em
 \vbox{\hsize3cm\tiny\raggedright\pretolerance10000
 \noindent #1\hfill}\hss}\vbox to8pt{\vfil}\vss}}}%
                        
\newcommand{\bgw}{{\textstyle \bigwedge}}
\newcommand{\bgs}{{\textstyle \bigodot}}
\newcommand{\bgt}{{\textstyle \bigotimes}}
\newcommand{\Span}{\mathrm{Span}}

\numberwithin{equation}{section}

\newcounter{romenumi}
\newcommand{\labelromenumi}{(\roman{romenumi})}

\newcommand{\ten}{\Upsilon}

\newcommand{\bma}{\begin{pmatrix}}
\newcommand{\ema}{\end{pmatrix}}
\newcommand{\gp}{{\stackrel{\scriptscriptstyle{+}}{\Gamma}}\phantom{}}
\newcommand{\gm}{{\stackrel{\scriptscriptstyle{-}}{\Gamma}}\phantom{}}
\newcommand{\ssp}{{\stackrel{\scriptscriptstyle{+}}{\Gamma}}_\spina\phantom{}}
\newcommand{\ssm}{{\stackrel{\scriptscriptstyle{-}}{\Gamma}}_\spina\phantom{}}
\newcommand{\kp}{{\stackrel{\scriptscriptstyle{+}}{\Omega}}\phantom{}}
\newcommand{\km}{{\stackrel{\scriptscriptstyle{-}}{\Omega}}\phantom{}}

\newcommand{\Rp}{{\stackrel{\scriptscriptstyle{+}}{R}}\phantom{}}
\newcommand{\Rm}{{\stackrel{\scriptscriptstyle{-}}{R}}\phantom{}}

\begin{document}
\title{Analog of selfduality in dimension nine} 
\vskip 1.truecm \author{Anna Fino} \address{Dipartimento di Matematica, Universit\`a di Torino, Via Carlo Alberto 10, 10123 Torino, Italy}
\email{annamaria.fino@unito.it}
\author{Pawe\l~ Nurowski} \address{Instytut Fizyki Teoretycznej,
Uniwersytet Warszawski, ul. Ho\.za 69, Warszawa, Poland}
\email{nurowski@fuw.edu.pl} \thanks{This reaserch was supported by the
  Polish Ministry of Research and Higher Education under grants NN201 607540 and NN202 104838, by  the MIUR project Differential Geometry and Global
Analysis (PRIN07) and by GNSAGA (Indam)}
\date{\today}

\begin{abstract}
We   introduce a type of Riemannian geometry in nine  dimensions,
which can be viewed as the counterpart of selfduality in   four
dimensions.  This geometry is   related to a $9$-dimensional
irreducible  representation of  ${\bf SO}(3) \times {\bf SO} (3)$ and
it turns out to be defined by  a  differential $4$-form. Structures
admitting a metric connection with totally antisymmetric  torsion
and preserving the $4$-form are studied in detail, producing  locally homogeneous examples which can be viewed as analogs of   self-dual  $4$-manifolds in dimension nine. 
\end{abstract}
\maketitle
\tableofcontents
\newcommand{\gat}{\tilde{\gamma}}
\newcommand{\Gat}{\tilde{\Gamma}}
\newcommand{\thet}{\tilde{\theta}}
\newcommand{\Thet}{\tilde{T}}
\newcommand{\rt}{\tilde{r}}
\newcommand{\st}{\sqrt{3}}
\newcommand{\kat}{\tilde{\kappa}}
\newcommand{\kz}{{K^{{~}^{\hskip-3.1mm\circ}}}}


\section{Introduction}

The special feature of $4$ dimensions is that the the rotation group  ${\bf SO} (4)$ is not simple but   it is locally isomorphic to ${\bf SU}(2) \times {\bf SU} (2)$, since 
${\mathfrak {so}}(4) = {\mathfrak {su}}(2)_L \oplus  {\mathfrak {su}}(2)_R$.

Given  an oriented $4$-dimensional Riemannian manifold $(M^4, g)$, the
Hodge-star-operator $* :\Lambda^2 \rightarrow  \Lambda^2$  satisfies
$*^2 = {\rm id}$ and the bundle  of $2$-forms $\Lambda^2$ splits as:
\begin{equation} \label{dec2-formsdim4}
\Lambda^2 = \Lambda^2_+ \oplus \Lambda^2_-,
\end{equation}
where $\Lambda^2_+$ is  the space of self-dual  forms and $\Lambda^2_-$ is  the  one of anti-self-dual  forms. 

The Riemann curvature tensor defines a self-adjoint transformation ${\mathcal R}: \Lambda^2 \rightarrow \Lambda^2$
%
which can be written, with respect  to the decomposition \eqref{dec2-formsdim4},  as  the  block matrix 
$$
{\mathcal R} = \left ( \begin{array}{cc} A&B\\
B^*&C \end{array} \right ),
$$
where $B \in {\mbox {Hom}}  (\Lambda^2_-, \Lambda^2_+)$ and $A \in  {\mbox {End}}  \, \Lambda^2_+$, $ C \in  {\mbox {End}}  \, \Lambda^2_-$ are self-adjoint.

This decomposition of ${\mathcal R}$  gives the  complete description of the Riemannian curvature tensor into irreducible components  obtained in  \cite{ST}:
$$
\left({\mbox {tr}} \, A, B,  A - \frac 13 {\mbox {tr}} \, A, C -  \frac 13 {\mbox {tr}} \, C \right),
$$
where ${\mbox {tr}} \, A = {\mbox {tr}} \, C$ is  the Ricci scalar, $B$ is the  traceless Ricci tensor, and the last  two components $W_+= A - \frac 13 {\mbox {tr}}  \, A$ and $W_- = C -  \frac 13 {\mbox {tr}} \, C$, together give the conformally invariant Weyl tensor $W = W_+ + W_-$.
We recall by   \cite{AHS}
 that  $g$ is Einstein if and only if $B = 0$ and $g$ is self-dual   if and only if $W_- = 0$.

In terms of   Lie algebra valued 1-form $\lc$   of the Levi-Civita connetion and of  its curvature  $2$-form ${{\stackrel{\scriptscriptstyle{LC}}{\Omega}}\phantom{}}$  we have the decompositions:
$$
\lc = \gp +  \gm, \quad \quad 
{{\stackrel{\scriptscriptstyle{LC}}{\Omega}}\phantom{}} =  {{\stackrel{\scriptscriptstyle{+}}{\Omega}}\phantom{}}  + {{\stackrel{\scriptscriptstyle{-}}{\Omega}}\phantom{}} ,
$$
where $\gp$ and $\kp$ are $\sua(2)_L$-valued, and $\gm$ and $\km$ are
$\sua(2)_R$-valued.

Then the condition for the Riemannian metric  $g$ to be   Einstein and self-dual   is equivalent to ${{\stackrel{\scriptscriptstyle{-}}{\Omega}}\phantom{}} =0$.

A natural problem  is to  study a  geometry in higher dimensions, which can be viewed as the counterpart of selfduality in   four dimensions.  The Lie group ${\bf SO}(n)$ for $n \geq 5$ is simple and there is no splitting of ${\mathfrak {so}}(n)$, so an idea is  to try  with a Lie group  of the form $H \times H$ in ${\bf SO}(n)$.

In this paper we will consider the case of ${\bf SO} (3) \times {\bf SO} (3) \subset  {\bf SO}(9)$. To this aim we need an irreducible $9$-dimensional representation of ${\bf SO} (3) \times {\bf SO} (3)$, which  turns out to be  related to   a  $9$-dimensional   irreducible representation  $\rho$  of the Lie group  $\slg(2,\bbR)\times\slg(2,\bbR)$. Perhaps for the first time 
the representation $\rho$  was used by G. Peano \cite{peano} in his extension of the
classical invariant theory to the action of the Cartesian product 
$\slg(2,\bbR)\times\slg(2,\bbR)$ on the Cartesian product $\bbR^2\times
\bbR^2$. Similarly to the classical invariant theory \cite[Ch. 10, p. 242]{olver}, Peano in
\cite{peano} defines irreducible representations of
$\slg(2,\bbR)\times\slg(2,\bbR)$ group, by considering its action on
homogeneous polynomials in  four  variables
$(\phi_1,\phi_2,\psi_1,\psi_2)=(\vec{\phi},\vec{\psi})\in\bbR^2\times\bbR^2$. Given a
defining action of $\slg(2,\bbR)$ on $\bbR^2$, $(h,\vec{\phi})\to
h\vec{\phi}$, the irreducible 
action of $\slg(2,\bbR)\times\slg(2,\bbR)$ on
$\bbR^{m+1}\times\bbR^{\mu+1}$, is defined as follows.

Let $a_{l\la}$, $l=0,\dots,m$, $\la=0,\dots,\mu$, be coordinates in
$\bbR^{m+1}\times\bbR^{\mu+1}$. They define 
a homogeneous polynomial
\be w(\vec{\phi},\vec{\psi})=\sum_{l=0}^m\sum_{\la=0}^\mu
a_{l\la}\binom{m}{l}\binom{\mu}{\la}\phi_1^{m-l}\phi_2^l\psi_1^{\mu-\la}\psi_2^\la.\label{poli}\ee
Now given $(h_L,h_R)\in \slg(2,\bbR)\times\slg(2,\bbR)$, we define
$a^{(h_L,h_R)}_{l\la}\in\bbR^{m+1}\times\bbR^{\mu+1}$ via:
$$\sum_{l=0}^m\sum_{\la=0}^\mu
a^{(h_L,h_R)}_{l\la}\binom{m}{l}\binom{\mu}{\la}\phi_1^{m-l}\phi_2^l\psi_1^{\mu-\la}\psi_2^\la=w(h_L\vec{\phi},h_R\vec{\psi}).$$
It follows that the map 
$$\slg(2,\bbR)\times\slg(2,\bbR)\times\bbR^{(m+1)(\mu+1)}\ni(h_L,h_R,a_{l\la})\to
(a^{(h_L,h_)}_{l\la})\in\bbR^{(m+1)(\mu+1)}$$
is an action of $\slg(2,\bbR)\times\slg(2,\bbR)$ on
$\bbR^{(m+1)(\mu+1)}$, and therefore it  
defines an $(m+1)(\mu+1)$-dimensional 
representation $\rho$ of this group by:
$$\rho(h_L,h_R)a_{l\la}=a^{(h_L,h_R)}_{l\la}.$$
For each value of $(m,\mu)$ this representation is irreducible.
In the paper  we are interested in the case $(m,\mu)=(2,2)$. In such
case the polynomial $w$ reads:
\be
\begin{array}{lcl}
w(\vec{\phi},\vec{\psi})&=&a_{00}\phi_1^2\psi_1^2+2a_{10}\phi_1\phi_2\psi_1^2+a_{20}\phi_2^2\psi_1^2+2a_{01}\phi_1^2\psi_1\psi_2+4a_{11}\phi_1\phi_2\psi_1\psi_2+\\[3 pt]
&&2a_{21}\phi_2^2\psi_1\psi_2+
a_{02}\phi_1^2\psi_2^2+2a_{12}\phi_1\phi_2\psi_2^2+a_{22}\phi_2^2\psi_2^2.
\end{array}\label{biq}
\ee
The 9-dimensional space $\bbR^9$ consisting of vectors 
$$\vec{x}=(x_0,x_1,x_2,x_3,x_4,x_5,x_6,x_7,x_8)=(a_{00},a_{10},a_{20},a_{01},a_{11},a_{21},a_{02},a_{12},a_{22}),$$
is equipped with the irreducible representation $\rho$ of
$\slg(2,\bbR)\times\slg(2,\bbR)$. This representation induces the
action of $\slg(2,\bbR)\times\slg(2,\bbR)$ on homogeneous polynomials
in variables $x_i$. Peano showed that the lowest order
\emph{invariant} polynomials under this action are:
\be
\begin{array}{c}
g=\displaystyle\sum_{i,j}g_{ij}x_ix_j=2\Big(x_0x_8+x_2x_6-2x_1x_7-2x_3x_5+2x_4^2\Big)\\[4 pt] \label{metold}
\ten=\displaystyle\sum_{i,j,k}\ten_{ijk}x_ix_jx_k =24\Big(x_0x_4x_8-x_0x_5x_7-x_1x_3x_8+x_1x_5x_6+
x_2x_3x_7-x_2x_4x_6\Big).
\end{array}
\ee

They equipp $\bbR^9$ with a metric $g_{ij}$ of signature $(4,5)$ and a
totally symmetric third rank tensor $\ten_{ijk}$, which turns out to be
traceless, $g^{ij}\ten_{ijk}=0$. 

The  common stabilizer of  the two  tensors $g$ and $\ten$, defined above, 
 is $\slg(2,\bbR)\times\slg(2,\bbR)$ in the
9-dimensional irreducible representation $\rho$ of Peano.

This is very similar to the situation 
in $\bbR^5$, where we have a pair of tensors $(g_{ij},\ten_{ijk})$
which reduce the $\glg(5,\bbR)$ group to the irreducible $\sog(3)$ in
dimension five \cite{ABF,bobi,cf}. The only difference with the 5-dimensional
case considered in \cite{bobi} is that there the metric $g_{ij}$ is of
\emph{purely Riemannian} signature \footnote{This indicates that the
  geometry associated with tensors $g$ and $\ten$ as above can be
  related to the geometry of a certain type of systems of differential equations of
  finite type \cite{god,nh}. Actually, the biquadrics (\ref{biq}) are
  related to the general solution of the finite type system $z_{xxx}=0$ $\&$
  $z_{yyy}=0$ of PDEs on the plane for the unknown $z=z(x,y)$. We
  expect that the geometry associated with $g$ and $\ten$ is the
  geometry of generalizations of this system \cite{dub}.}; see also \cite{gut,nur0,nurm1}. 

The Riemannian version of tensors associated with Peano biquadrics may 
be obatined by making
the following formal substitutions in \eqref{metold}:
$$
\begin{array}{lll}
x_0=y_1+i y_2, &x_8=y_1-i y_2, & x_2=y_3+i y_4\\[3pt]
x_6=y_3-iy_4, &x_1=\tfrac{1}{\sqrt{2}}(y_5+iy_6),
& x_7=-\tfrac{1}{\sqrt{2}}(y_5-iy_6)\\[3 pt]
x_3=\tfrac{1}{\sqrt{2}}(y_7+iy_8), &
x_5=-\tfrac{1}{\sqrt{2}}(y_7-iy_8), 
& x_4=\tfrac{1}{\sqrt{2}}y_9.
\end{array}  
$$
In these formulae coefficients $y_\mu$, $\mu=1,\dots,9$, are
\emph{real}, and $i$ is the 
imaginary unit. With these substitutions  \eqref{metold}
become:
\begin{equation}\label{newmet}
\begin{array}{lcl}
g &=&\displaystyle\sum_{i,j}g_{ij}y_iy_j=2\Big(y_1^2+y_2^2+y_3^2+y_4^2+y_5^2+y_6^2+y_7^2+y_8^2+y_9^2\Big),\\[4 pt]
\ten &= &\displaystyle\sum_{i,j,k}\ten_{ijk}y_iy_jy_k=\\[4 pt]
&&12\Big(-2y_1y_5y_7-2y_3y_5y_7-2y_2y_6y_7-2y_4y_6y_7-2y_2y_5y_8+\\
&&2y_4y_5y_8+2y_1y_6y_8-2y_3y_6y_8+\sqrt{2}y_1^2y_9+\sqrt{2}y_2^2y_9-\sqrt{2}y_3^2y_9-\sqrt{2}y_4^2y_9\Big).
\end{array}
\end{equation}
This equipps $\bbR^9$ parametrized by $y_\mu$, $\mu=1, 2, \dots,9,$ with a
pair of totally symmetric tensors $(g_{ij},\ten_{ijk})$, in which
$g_{ij}$ is now a Riemannian  metric. 

In Section \ref{sectioninvariant}  we obtain a better realization of 
$(\bbR^9,  g, \ten)$ by using the indentification of $\bbR^9$ with  a space
$\bbM_{3\times 3}(\bbR)$ of $3\times3$ matrices with real coefficients. This allows us to show that 
${\bf SO}(3)\times{\bf SO}(3)$ is   surprising the stabilizer of  a $4$-form $\omega$.  In Section \ref{repr} irreducible representations of ${\bf SO}(3)\times{\bf SO}(3)$ are studied in detail.   Following  the approach presented in \cite{bobi}, in Section \ref{secirrgeom} we introduce the irreducible $\sog(3)\times\sog(3)$ geometry in dimension nine as 
the geometry of 
9-dimensional manifolds $M^9$ equipped either  with  a  pair of totally symmetric
tensors $(g,\ten)$ as in \eqref{newmet} or with  the differential  $4$-form $\omega$.  In Section \ref{secnearlyint} we 
 determine the conditions for $\ten$  which will
guarantee that  $(M^9,g,\ten, \omega)$ admits a unique metric
connection $\Gamma$, with values in the symmetry algebra $(\soa(3)_L\oplus\soa(3)_R)$ of  $(g,\ten)$ and with  totally antisymmetric torsion. 
This  $(\soa(3)_L\oplus\soa(3)_R)$-connection $\Gamma$, also called the  characteristic connection,  naturally splits onto: 
$$\Gamma=\gp+\gm, 
$$
  with $\gp\in\soa(3)_L\otimes\R^9,$
and $\gm\in\soa(3)_R\otimes\R^9.$ Because $\soa(3)_L$ commutes with $\soa(3)_R$ this split defines two
independent $\soa(3)$-valued connections $\gp$ and $\gm$. So an irreducible $\sog(3)\times\sog(3)$ geometry
$(M^9,g,\ten,\omega)$ equipped with and $(\soa(3)_L\oplus\soa(3)_R)$
connection $\Gamma$ can be Einstein in several meanings, by considering not only the Levi-Civita connection but also the connections $\Gamma$, $\gp$ and $\gm$. In the last section  we study   irreducible $\sog(3)\times\sog(3)$ geometries
$(M^9,g,\ten,\omega)$   admitting a characteristic  connection $\Gamma$ with   `special'  torsion  $T$. In particular,  we provide locally homogeneous  (non Riemannian symmetric) examples   for which $T \neq 0$,  $\gp$ has  vanishing curvature  and  $\gm$ is Einstein and not flat.  These examples can be viewed as analogs of self-dual structures in dimension four. It would be very interesting to find
analogs of selfduality which are not locally homogeneous. If such
solutions may exist is an open question. 
 
\medskip

\noindent {\it Acknowledgments.}
We would like to thank Robert Bryant,  Antonio Di Scala, Boris Doubrov, Mike Eastwood,  Katja Sagerschnig and  Simon Salamon  for useful comments and suggestions. 

\medskip



\section{Invariant ${\bf SO}(3)\times
  {\bf SO}(3)$ tensors} \label{sectioninvariant} 
We identify the 9-dimensional real vector space $\bbR^9$ with a space
$\bbM_{3\times 3}(\bbR)$ of $3\times3$ matrices with real coefficients,  via the map $$\sigma:\bbR^9\to \bbM_{3\times 3}(\bbR),$$ defined by
\be
\bbR^9\ni A=a^i\bas_i\longmapsto \sigma(A)=\begin{pmatrix} 
a^1&a^2&a^3\\
a^4&a^5&a^6\\
a^7&a^8&a^9
\end{pmatrix}\in \bbM_{3\times3}(\bbR).\label{bass}
\ee
This map is obviously invertible, so we also have the inverse 
$$\sigma^{-1}: \bbM_{3\times3}(\bbR)\to\bbR^9.$$
The unique irreducible 9-dimensional representation
$\rho$ of the group 
$$G={\bf SO}(3)\times{\bf SO}(3)$$ in $\mathbb{R}^9$ is then defined
as follows. 

Let $h=(h_L,h_R)$ be the most general element of $G$,
i.e. let $h_L$ and $h_R$ be two arbitrary elements of
$\sog(3)$ in the standard representation of $3\times 3$ real matrices. Then, 
for every vector $A$ from $\bbR^9$, we have: 
\be
\rho(h)A~=~\sigma^{-1}\big(h_L~\sigma(A)~h_R^{-1}\big).\label{so35}
\ee

In the rest of the article we adopt the convention that the symbol $G$
is reserved to denote the group $\sog(3)\times\sog(3)$ in the
irreducible 9-dimensional representation defined above, and that
$\mathfrak{g}$ denotes its Lie algebra, ${\mathfrak g}=\soa(3)\times\soa(3)$. 

Consider now $\theta=(\theta^1,\theta^2,\dots,\theta^9)$
with components $\theta^i$ being covectors in $\bbR^9$. This means
that $\theta$ is a vector-valued 1-form,
$\theta\in\bbR^9\otimes(\bbR^9)^*$.  We identify it
with the matrix-valued 1-form $$\sigma(\theta)\in\bbM_{3\times 3}((\bbR^9)^*).$$ The group
$G$ acts on forms $\theta$ via
$$\theta\mapsto \theta'=\rho(h)\theta.$$
Its action is then extended to all tensors $T$ of the form 
$$T=T_{i_1i_2...i_r}\theta^{i_1}\otimes\theta^{i_2}\otimes ...\otimes\theta^{i_r}$$
via 
$$T\mapsto
T'=T_{i_1i_2...i_r}{\theta'}^{i_1}\otimes{\theta'}^{i_2}\otimes ...\otimes{\theta'}^{i_r}.$$
We say that the tensor $T$ is $G$-invariant iff
$T'=T.$

An example of a $G$-invariant tensor is obtained by considering the 
determinant 
$${\rm det}(\sigma(A))=\tfrac16\ten_{ijk}a^ia^ja^k$$
and its corresponding \emph{symmetric} tensor
\be\ten:=\tfrac16\ten_{ijk}\theta^i\odot\theta^j\odot\theta^k.\label{de}\ee
This is obviously $G$-invariant by the properties of the determinant,
and by the fact that $\det(h)=1$, for every element of $\sog(3)$. 

Thus we have at least one $G$-invariant tensor $\ten$.

To create
others we note the $G$-invariance of the expressions  
\be
\hbox{Tr}(\sigma(\theta)\odot\sigma(\theta)^T),\quad\quad \hbox{Tr}(\sigma(\theta)\dz \sigma(\theta)^T),\quad\quad \hbox{Tr}(\sigma(\theta)\otimes \sigma(\theta)^T).\label{ek}\ee 
Here, the product sign under the trace is
considered as the usual row-by-columns product of
$3\times 3$ matrices, but with the product between the matrix
elements in each sum being the respective tensor
products $\odot$, $\dz$ and $\otimes$. 
The $G$-invariance of these three expressions is an immediate
consequence of the defining property of the elements of $\sog(3)$,
namely: $h^Th=hh^T=\id$. 
Having observed this, we now see that any function $F$, multilinear in
expressions (\ref{ek}), also defines a $G$-invariant tensor.

This enables us to
define a new $\sog(3)\times\sog(3)$-invariant tensor:
\be
g=\hbox{Tr}(\sigma(\theta)\odot\sigma(\theta)^T)=g_{ij}\theta^i\theta^j.\label{met}\ee
This tensor is symmetric, rank $\binom{0}{2}$ and nondegenerate. It
defines a Riemannian metric $g$ on $\bbR^9$. 

Another set of $G$-invariant tensors is given by the $2k$-forms
\be\hbox{Tr}(\sigma(\theta)\dz \sigma(\theta)^T\dz \sigma(\theta)\dz \sigma(\theta)^T\dz...\dz \sigma(\theta)\dz \sigma(\theta)^T).\label{ff1}\ee
One would expect that these identically vanish, but surprisingly, we have the
following proposition.
\begin{proposition}
The 4-form
\be
\omega=\tfrac14{\rm Tr}(\sigma(\theta)\dz \sigma(\theta)^T\dz \sigma(\theta)\dz \sigma(\theta)^T)=\tfrac{1}{4!}\omega_{ijkl}\theta^i\dz\theta^j\dz\theta^k\dz\theta^l\label{omeg}\ee
does \emph{not} vanish, $\omega\neq 0$. 

In the remaining cases, when $k=1$, 3, 4, the forms (\ref{ff1}) are
identically equal to zero.
\end{proposition}
We have the following formulae for the three $G$-invariant objects
defined above:
\be\begin{aligned}
\ten&=-\theta^3\theta^5\theta^7+\theta^2\theta^6\theta^7+\theta^3\theta^4\theta^8-\theta^1\theta^6\theta^8-\theta^2\theta^4\theta^9+\theta^1\theta^5\theta^9,\\
\\
g&=(\theta^1)^2+(\theta^2)^2+(\theta^3)^2+(\theta^4)^2+(\theta^5)^2+(\theta^6)^2+(\theta^7)^2+(\theta^8)^2+(\theta^9)^2,\\
\\
\omega&=\theta^1\dz\theta^2\dz\theta^4\dz\theta^5+\theta^1\dz\theta^2\dz\theta^7\dz\theta^8+\theta^1\dz\theta^3\dz\theta^4\dz\theta^6+\\&\theta^1\dz\theta^3\dz\theta^7\dz\theta^9+\theta^2\dz\theta^3\dz\theta^5\dz\theta^6+\theta^2\dz\theta^3\dz\theta^8\dz\theta^9+\\&\theta^4\dz\theta^5\dz\theta^7\dz\theta^8+\theta^4\dz\theta^6\dz\theta^7\dz\theta^9+\theta^5\dz\theta^6\dz\theta^8\dz\theta^9.\end{aligned}\label{ty}\ee
Here, to simplify the notation, we abreviated expressions like 
$\theta^3\odot\theta^5\odot\theta^7$, or $\theta^1\odot\theta^1$, to
the respective, $\theta^3\theta^5\theta^7$ and $(\theta^1)^2$.

We have the following   proposition.
\begin{proposition}$\phantom{}$\newline\vspace{-0.5cm}
\begin{enumerate}

\item The simultaneous stabilizer in $\glg(9,\bbR)$ of the tensors
$g$ and $\ten$ defined respectively in (\ref{de}) and (\ref{met}) is
$G=\sog(3)\times\sog(3)$ in the irreducible 9-dimensional
representation $\rho$.

 \item The stabilizer in $\glg(9,\bbR)$ of the 4-form $\omega$ defined in
 (\ref{omeg}) is also 
$G=\sog(3)\times\sog(3)$ in the irreducible 9-dimensional
representation $\rho$.
\end{enumerate}

\end{proposition}

\begin{proof} We know from the considerations preceeding the proposition that the
stabilizers contain $G$. To show that they are actually equal to $G$
we do as follows:

A stabilizer $G'$ of $g$ and $\ten$ consists of those 
elements $h$ in $\glg(9,R)$ for which
\be
g(hX,hY)=g(X,Y)\quad\quad{\rm and}\quad\quad \ten(hX,hY,hZ)=\ten(X,Y,Z).\label{re}\ee  
We find the Lie algebra of $G'$. Taking $h$ in the form $h=\exp(sX)$
and taking $\tfrac{\rm d}{{\rm d}s}_{|s=0}$  of the equations
(\ref{re}), we see that the matrices $X=(X^i_{~j})$ representing the
elements of the Lie algebra $\mathfrak{g}'$ of $G'$ must satisfy
\be
g_{lj}X^{l}_{~~i}+g_{il}X^{l}_{~~j}=0\label{rg}
\ee
and
\be
\ten_{ljk}X^{l}_{~~i}+\ten_{ilk}X^{l}_{~~j}+\ten_{ijl}X^{l}_{~~k}=0\label{rt}.
\ee
The first of the above equations tells that the matrices $X$ must be
antisymmetric, i.e. it reduces 81 components of a matrix $X$ to 36. The
second equation gives another 30 independent conditions restricting
the number of free components of $X$ to 6. Explicitly the matrix $X$
solving (\ref{rg})-(\ref{rt}) is of the form 
\be X=X^1e_1+X^2e_2+X^3e_3+X^{1'}e_{1'}+X^{2'}e_{2'}+X^{3'}e_{3'},\label{lieg}\ee
where
\be
\begin{aligned}
&e_1=\left(\begin{smallmatrix} 
0&0&0&-1&0&0&0&0&0\\
0&0&0&0&-1&0&0&0&0\\
0&0&0&0&0&-1&0&0&0\\
1&0&0&0&0&0&0&0&0\\
0&1&0&0&0&0&0&0&0\\
0&0&1&0&0&0&0&0&0\\
0&0&0&0&0&0&0&0&0\\
0&0&0&0&0&0&0&0&0\\
0&0&0&0&0&0&0&0&0\end{smallmatrix}\right),e_2=\left(\begin{smallmatrix} 
0&0&0&0&0&0&-1&0&0\\
0&0&0&0&0&0&0&-1&0\\
0&0&0&0&0&0&0&0&-1\\
0&0&0&0&0&0&0&0&0\\
0&0&0&0&0&0&0&0&0\\
0&0&0&0&0&0&0&0&0\\
1&0&0&0&0&0&0&0&0\\
0&1&0&0&0&0&0&0&0\\
0&0&1&0&0&0&0&0&0\end{smallmatrix}\right)
,e_3=\left(\begin{smallmatrix} 
0&0&0&0&0&0&0&0&0\\
0&0&0&0&0&0&0&0&0\\
0&0&0&0&0&0&0&0&0\\
0&0&0&0&0&0&-1&0&0\\
0&0&0&0&0&0&0&-1&0\\
0&0&0&0&0&0&0&0&-1\\
0&0&0&1&0&0&0&0&0\\
0&0&0&0&1&0&0&0&0\\
0&0&0&0&0&1&0&0&0\end{smallmatrix}\right),
\\
& e_{1'}=\left(\begin{smallmatrix} 
0&-1&0&0&0&0&0&0&0\\
1&0&0&0&0&0&0&0&0\\
0&0&0&0&0&0&0&0&0\\
0&0&0&0&-1&0&0&0&0\\
0&0&0&1&0&0&0&0&0\\
0&0&0&0&0&0&0&0&0\\
0&0&0&0&0&0&0&-1&0\\
0&0&0&0&0&0&1&0&0\\
0&0&0&0&0&0&0&0&0\end{smallmatrix}\right),
e_{2'}=\left(\begin{smallmatrix} 
0&0&-1&0&0&0&0&0&0\\
0&0&0&0&0&0&0&0&0\\
1&0&0&0&0&0&0&0&0\\
0&0&0&0&0&-1&0&0&0\\
0&0&0&0&0&0&0&0&0\\
0&0&0&1&0&0&0&0&0\\
0&0&0&0&0&0&0&0&-1\\
0&0&0&0&0&0&0&0&0\\
0&0&0&0&0&0&1&0&0\end{smallmatrix}\right),
e_{3'}=\left(\begin{smallmatrix} 
0&0&0&0&0&0&0&0&0\\
0&0&-1&0&0&0&0&0&0\\
0&1&0&0&0&0&0&0&0\\
0&0&0&0&0&0&0&0&0\\
0&0&0&0&0&-1&0&0&0\\
0&0&0&0&1&0&0&0&0\\
0&0&0&0&0&0&0&0&0\\
0&0&0&0&0&0&0&0&-1\\
0&0&0&0&0&0&0&1&0\end{smallmatrix}\right).\end{aligned}\label{dwie}\ee
It is easy to check that the matrices $e$ satisfy the following
commutation relations: $[e_1,e_2]=e_3$, $[e_3,e_1]=e_2$,
$[e_2,e_3]=e_1$, $[e_{1'},e_{2'}]=e_{3'}$, $[e_{3'},e_{1'}]=e_{2'}$,
$[e_{2'},e_{3'}]=e_{1'}$, with all the other commutators being zero
modulo the antisymmetry. Thus the system $(e_A,e_{A'})$, $A=1,2,3$,
spans the Lie algebra  $\soa(3)\oplus\soa(3)$, confirming that 
the Lie algebra
$\mathfrak{g}'$ of the stabilizer $G'$ of tensors (\ref{de}) and
(\ref{met}) is $\mathfrak{g}'=\soa(3)\oplus\soa(3)$. 
In an analogous way we find the Lie algebra $\mathfrak{g}''$ of the
stabilizer $G''$ of $\om$. 
This stabilizer consists of those 
elements $h$ in $\glg(9,\bbR)$ for which
\be
\om(hX,hY,hZ)=\om(X,Y,Z).\label{rep}\ee  
Taking $h$ in the form $h=\exp(sX)$
and taking $\tfrac{\rm d}{{\rm d}s}_{|s=0}$  of the equations
(\ref{rep}), we see that the matrices $X=(X^i_{~j})$ representing the
elements of the Lie algebra $\mathfrak{g}''$ of $G''$ must satisfy
\be
\omega_{ljkm}X^{l}_{~~i}+\omega_{ilkm}X^{l}_{~~j}+\omega_{ijlm}X^{l}_{~~k}+\omega_{ijkl}X^{l}_{~~m}=0\label{ro}.
\ee
A short algebra shows that this imposes 75 independent conditions on
the 81 components of $X$, and that the most general solution to this
equation is given by (\ref{lieg}) with the generators $(e_A,e_{A'})$
as in (\ref{dwie}). Thus 
$$\mathfrak{g}'=\mathfrak{g}''=\soa(3)\oplus\soa(3):=\mathfrak{g}.$$
As a consequence $G' = G'' = \sog(3)\times\sog(3)$,  since $\soa(3)\oplus\soa(3)$ is  a maximal Lie subalgebra of $\soa(9)$.
\end{proof}
\begin{remark}
Note that the form $\omega$ \emph{alone} is enough to reduce
$\glg(9,\bbR)$ to $G$. One does not need the metric $g$ for this reduction! 
On the other hand,   the tensor $\ten$ alone is not enough to reduce the
$\glg(9,\bbR)$ to $G$. The equation (\ref{rt}) imposes only 65
independent conditions on the matrix $X$. Thus it reduces
$\gla(9,\bbR)$ to a Lie algebra of dimension 16. Since 16 is the
dimension of $\sla(3,\bbR)\oplus\sla(3,\bbR)$, and $\ten$ is clearly
  $\slg(3,\bbR)\times\slg(3,\bbR)$-invariant, the stabilizer of the
  tensor $\ten$ alone is $\slg(3,\bbR)\times\slg(3,\bbR)$. To reduce
  it further to
$\soa(3)\oplus\soa(3)$ one needs to preserve $g$. If in addition to
$\ten$ we preserve $g$ we get, via the equation (\ref{rg}), the remaining 10 conditions.
\end{remark}
\begin{remark} For the geometric relevance of the form $\omega$ see Remark \ref{Pontriagyn} suggested by Robert Bryant  (\cite{br,mic}).
\end{remark}
\begin{remark}
We remark that in addition to the 4-form $\omega$ we have also the
5-form $*\omega$  (Hodge-dual of $\omega$)  which is $G$-invariant.  One can say that given only $\omega$ in
$\bbR^9$ we do not have any metric structure on it. But $\omega$
defines the reduction of the Lie algebra of $\glg(9,\bbR)$ to
$\mathfrak{g}=\soa(3)\times\soa(3)$. In particular it defines the explicit
representation of $\mathfrak{g}$ given by (\ref{lieg}) with the
explicit form of the generators $(e_A,e_{A'})$ given by (\ref{dwie}). Thus, given $\omega$
we have explicitly $X$ as in (\ref{lieg}). 
Now we define the
metric $g_{ij}$ as a $\binom{0}{2}$-tensor such
that (\ref{omeg}) holds. It is a matter of checking that given
$X$ as in (\ref{lieg}) with $(e_A,e_{A'})$ as in (\ref{dwie}) the only metric
$g_{ij}$ satisfying (\ref{omeg}) (miraculosly!) is $g_{ij}={\rm
  const}\times \delta_{ij}$. Thus the 4-form $\omega$ defines the metric $g$ up to
a scale, and this in turn defines the unique (up to a scale) 5-form
$*\omega$, being its standard Hodge-star with respect to the metric
$g$.   
\end{remark}
Another way of defining the 5-form $*\omega$, which provides the
explicit relation between $(g,\ten)$ and $\omega$, is given by the
proposition below. To formulate it we consider a coframe $\theta^i$
and the corresponding components $\ten_{ijk}$ of the tensor
$\ten$ as in (\ref{de}). Using them we define a $(9\times 9)$-matrix-valued-1-form
$\ten(\theta)=(\ten(\theta)^i_{~j})$ with matrix elements 
$$\ten(\theta)^i_{~j}=g^{il}\ten_{ljk}\theta^k.$$
Here $(g^{ij})$ is the matrix inverse of $(g_{ij})$,
i.e. $g^{ik}g_{kj}=g_{jk}g^{ki}=\delta^i_{~j}$.  
Having the matrix $\ten(\theta)$, we consider traces of the skew symmetric powers of it,
$${\rm Tr}(\ten(\theta)^{\dz k})={\rm Tr}(\ten(\theta)\dz\ten(\theta)\dz...\dz\ten(\theta)),$$
where again the expressions like $\ten(\theta)\dz\ten(\theta)$ denote
the usual row-by-columns multiplication of $9\times 9$ matrices, with
the multiplication between the matrix elements being the wedge product
$\dz$.  
\begin{proposition}
If $k\neq 5$ and $k\in\{1,2,...,9\}$, then ${\rm Tr}(\ten(\theta)^{\dz k})=0$.

If $k=5$ the 5-form ${\rm Tr}(\ten(\theta)^{\dz 5})$ does not vanish,
$${\rm Tr}(\ten(\theta)^{\dz 5})={\rm
  Tr}(\ten(\theta)\dz\ten(\theta)\dz\ten(\theta)\dz\ten(\theta)\dz\ten(\theta))\neq
0.$$
Up to a scale this form is equal to the $G$-invariant 5-form $*\omega$.
In turn, the relation between the form $\omega$ and tensors $(g,\ten)$
is given by
$$\omega=*{\rm
  Tr}(\ten(\theta)\dz\ten(\theta)\dz\ten(\theta)\dz\ten(\theta)\dz\ten(\theta)).$$
\end{proposition} 
We proved this proposition by a brute force, using 
\eqref{ty}, and calculating  the expression of ${\rm
  Tr}(\ten(\theta)^{\dz k})$ for each value of $k=1,2,...9$. 
It would be interesting to get a `pure thought' proof of it.
\begin{remark}
The situation with $G$-invariant totally \emph{anti}symmetric $p$-forms is clear: there are only
one (up to a scale) 0- and 9-forms (a constant and its Hodge dual), and there are only one (up to a
scale) 4- and 5-forms (the 4-form $\omega$ and its Hodge dual). All the other
$G$-invariant $p$-forms are equal to zero.
\end{remark}
\begin{remark}
The situation with $G$-invariant totally \emph{symmetric} 
$p$-forms is more complex because of the infinite dimension of
$\bigoplus_{k=0}^\infty\bigodot^k\bbR^9$: Up to a scale there is only
one totally symmetric $G$-invariant 0-form; totally
symmetric $G$-invariant 1-forms are all equal to zero; there is only one totally symmetric $G$-invariant 2-form - the metric
$g$, and only one totally symmetric $G$-invariant 3-form - the tensor
$\ten$. Continuing this one gets that, in particular, there is only 
a 2-real-parameter family of totally symmetric $G$-invariant 4-forms: 
the family is spanned by $g_{(ij}g_{kl)}$ and by a 
tensor $\Xi_{ijkl}=\Xi_{(ijkl)}$, which in our coframe $\theta$ is
expressed by:
$$\begin{aligned}
\Xi=&\tfrac{1}{24}\Xi_{ijkl}\theta^i\theta^j\theta^k\theta^l=\\
&2 (\theta^1)^4+4 (\theta^1)^2 (\theta^2)^2+2 (\theta^2)^4+4 (\theta^1)^2 (\theta^3)^2+4 (\theta^2)^2 (\theta^3)^2+2 (\theta^3)^4+\\&4 (\theta^1)^2 (\theta^4)^2-7 (\theta^2)^2 (\theta^4)^2-7 (\theta^3)^2 (\theta^4)^2+2 (\theta^4)^4+22 \theta^1 \theta^2 \theta^4 \theta^5-\\&7 (\theta^1)^2 (\theta^5)^2+4 (\theta^2)^2 (\theta^5)^2-7 (\theta^3)^2 (\theta^5)^2+4 (\theta^4)^2 (\theta^5)^2+2 (\theta^5)^4+\\&22 \theta^1 \theta^3 \theta^4 \theta^6+22 \theta^2 \theta^3 \theta^5 \theta^6-7 (\theta^1)^2 (\theta^6)^2-7 (\theta^2)^2 (\theta^6)^2+4 (\theta^3)^2 (\theta^6)^2+\\&4 (\theta^4)^2 (\theta^6)^2+4 (\theta^5)^2 (\theta^6)^2+2 (\theta^6)^4+4 (\theta^1)^2 (\theta^7)^2-7 (\theta^2)^2 (\theta^7)^2-\\&7 (\theta^3)^2 (\theta^7)^2+4 (\theta^4)^2 (\theta^7)^2-7 (\theta^5)^2 (\theta^7)^2-7 (\theta^6)^2 (\theta^7)^2+2 (\theta^7)^4+\\&22 \theta^1 \theta^2 \theta^7 \theta^8+22 \theta^4 \theta^5 \theta^7 \theta^8-7 (\theta^1)^2 (\theta^8)^2+4 (\theta^2)^2 (\theta^8)^2-\\&7 (\theta^3)^2 (\theta^8)^2-7 (\theta^4)^2 (\theta^8)^2+4 (\theta^5)^2 (\theta^8)^2-7 (\theta^6)^2 (\theta^8)^2+4 (\theta^7)^2 (\theta^8)^2+\\&2 (\theta^8)^4+22 \theta^1 \theta^3 \theta^7 \theta^9+22 \theta^4 \theta^6 \theta^7 \theta^9+22 \theta^2 \theta^3 \theta^8 \theta^9+22 \theta^5 \theta^6 \theta^8 \theta^9-\\&7 (\theta^1)^2 (\theta^9)^2-7 (\theta^2)^2 (\theta^9)^2+4 (\theta^3)^2 (\theta^9)^2-7 (\theta^4)^2 (\theta^9)^2-7 (\theta^5)^2 (\theta^9)^2+\\&4 (\theta^6)^2 (\theta^9)^2+4 (\theta^7)^2 (\theta^9)^2+4 (\theta^8)^2 (\theta^9)^2+2 (\theta^9)^4.
\end{aligned}
$$
The $G$-invariant tensor $\Xi_{ijkl}$ defined above may be
characterized as the unique (up to a scale) $G$-invariant totally symmetric
$\binom{0}{4}$ tensor which has vanishing trace, $g^{ij}\Xi_{ijkl}=0$.   
\end{remark}
\section{Irreducible representations of $\sog(3)\times\sog(3)$}\label{repr}
As it is well known all finite dimensional real irreducible representations
of $\sog(3)$ have dimensions $d_k=2k+1$, $k=0,1,2,3,....$, and are
enumerated by the weight vectors $[2k]$. The representations with
the weight vectors $[m]=[2k]$ and $[\mu]=[2l]$ are
equivalent\footnote{Note that $m$ and $\mu$ here are related to the
  order of the Peano plynomials in (\ref{poli}).} iff
$k=l$. We
denote the vector spaces of these representations by $V_{[2k]}$.  
Consequently, all pairwaise inequivalent finite dimensional real
irreducible representations of $\sog(3)\times\sog(3)$ are given by
tensor products 
$$V_{[2k]}\otimes V_{[2l]}:=V_{[2k,2l]},\quad\quad{\rm with}\quad k,l=0,1,2,3,...,$$
and have the respective dimensions 
$$d_{[2k,2l]}=(2k+1)(2l+1).$$ 
In particular, for each number $d_{[2k,2l]}$, with $k\neq l$, there are two
nonequivalent irreducible representations of $\sog(3)\times\sog(3)$
with the respective carrier spaces $V_{[2k,2l]}$ and $V_{[2l,2k]}$.

In the following we will need decompositions of various tensor
products of spaces $V_{[2k,2l]}$ onto irreducible components with respect to the action of $\sog(3)\times\sog(3)$. These are summarized in

\begin{proposition}\label{dec1}
$$\bgw^2 V_{[2,2]}=V_{[0,2]}\oplus V_{[2,0]}\oplus V_{[2,4]}\oplus V_{[4,2]},$$
$$\bgw^3 V_{[2,2]}=V_{[0,2]}\oplus V_{[2,0]}\oplus V_{[0,6]}\oplus
  V_{[6,0]}\oplus V_{[2,4]}\oplus V_{[4,2]}\oplus V_{[2,2]}\oplus
  V_{[4,4]},$$
$$\bgw^4 V_{[2,2]}=V_{[0,0]}\oplus V_{[0,4]}\oplus V_{[4,0]}\oplus 2V_{[2,2]}\oplus
  V_{[2,4]}\oplus V_{[4,2]}\oplus V_{[2,6]}\oplus V_{[6,2]}\oplus
  V_{[4,4]},$$
$$\bgs^2 V_{[2,2]}=V_{[0,0]}\oplus V_{[0,4]}\oplus V_{[4,0]}\oplus
  V_{[2,2]}\oplus V_{[4,4]},$$
$$\bgs^3 V_{[2,2]}=V_{[0,0]}\oplus 2V_{[2,2]}\oplus V_{[2,4]}\oplus
  V_{[4,2]}\oplus V_{[2,6]}\oplus V_{[6,2]}\oplus V_{[4,4]}\oplus
  V_{[6,6]},$$
$$
\begin{aligned}
\bgs^4 V_{[2,2]}=2V_{[0,0]}\oplus 2V_{[0,4]}\oplus 2V_{[4,0]}\oplus
2V_{[2,2]}\oplus V_{[2,4]}\oplus  V_{[4,2]}\oplus V_{[0,8]}\oplus V_{[8,0]}\oplus\\
  V_{[2,6]}\oplus V_{[6,2]}\oplus 3 V_{[4,4]}\oplus V_{[4,6]}\oplus
  V_{[6,4]}\oplus V_{[4,8]}\oplus V_{[8,4]}\oplus V_{[6,6]}\oplus V_{[8,8]}.
\end{aligned}
$$
\end{proposition}
We in addition have the following identifications:
$$
\begin{aligned}
{\rm `left'}\quad\soa(3)=V_{[0,2]}\\
{\rm `right'}\quad\soa(3)=V_{[2,0]}\\
\bbR^9=V_{[2,2]}\\ 
\soa(9)=\bgw^2\bbR^9=\bgw^2V_{[2,2]}.
\end{aligned}
$$
In the following we will conveniently denote the $\soa(3)$ Lie
algebra corresponding to $V_{[0,2]}$ by $\soa(3)_L$ and the $\soa(3)$ Lie
algebra corresponding to $V_{[2,0]}$ by $\soa(3)_R$, i.e.
$$V_{[0,2]}=\soa(3)_L,\quad\quad{\rm and}\quad\quad V_{[2,0]}=\soa(3)_R.$$ 
Using these identifications and the decompositions from the
proposition above, we obtain:
\begin{proposition}~\\
$$
\begin{aligned}
\soa(9)\otimes \bbR^9&=2V_{[0,2]}\oplus 2V_{[2,0]}\oplus V_{[0,4]}\oplus
V_{[4,0]}\oplus  V_{[0,6]}\oplus V_{[6,0]}\oplus\\ &3V_{[2,4]}\oplus 3V_{[4,2]}\oplus V_{[2,6]}\oplus
  V_{[6,2]}\oplus 4V_{[2,2]}\oplus2V_{[4,4]}\oplus\\& V_{[4,6]}\oplus
  V_{[6,4]}\\
\soa(3)_L\otimes \bbR^9&=V_{[2,0]}\oplus V_{[2,2]}\oplus V_{[2,4]}\\
\soa(3)_R\otimes \bbR^9&=V_{[0,2]}\oplus V_{[2,2]}\oplus V_{[4,2]}\\
\Big(\soa(3)_L\oplus\soa(3)_R\Big)\otimes \bbR^9&=V_{[0,2]}\oplus
V_{[2,0]}\oplus 2V_{[2,2]}\oplus V_{[2,4]}\oplus V_{[4,2]}\\
\soa(3)_L\otimes \bgw^2\bbR^9&=V_{[0,0]}\oplus
V_{[0,2]}\oplus V_{[2,6]}\oplus V_{[0,4]}\oplus
V_{[4,0]}\oplus V_{[2,4]}\oplus V_{[4,2]}\oplus\\&
2V_{[2,2]}\oplus V_{[4,4]}\\
\soa(3)_R\otimes \bgw^2\bbR^9&=V_{[0,0]}\oplus
V_{[2,0]}\oplus V_{[6,2]}\oplus V_{[0,4]}\oplus
V_{[4,0]}\oplus V_{[2,4]}\oplus V_{[4,2]}\oplus\\& 2V_{[2,2]}\oplus V_{[4,4]}\\
\Big(\soa(3)_L\oplus\soa(3)_R\Big)\otimes \bgw^2\bbR^9&=2V_{[0,0]}\oplus
V_{[0,2]}\oplus V_{[2,0]}\oplus 2V_{[0,4]}\oplus
2V_{[4,0]}\oplus\\&2V_{[2,4]}\oplus2V_{[4,2]}\oplus V_{[2,6]}\oplus V_{[6,2]}\oplus 4V_{[2,2]}\oplus2V_{[4,4]}\end{aligned}
$$
\end{proposition}
The proofs of the above propositions can be
obtained by the standard representation theory methods using 
weights. Instead of presenting them we identify various useful
components of 
the decompositions mentioned in the propositions as 
eigenspaces of certain $\sog(3)\times\sog(3)$ invariant operators.  

For example the four irreducible components in the decomposition of
$\bgw^2\bbR^9$ in Proposition \ref{dec1} can be distinguished by means
of the action of the endomorphism of $\bgt^2\bbR^9$ defined by the
structural 4-form $\omega$. Indeed the 4-form
$\omega=\tfrac{1}{24}\omega_{ijkl}\theta^i\dz\theta^j\dz\theta^k\dz\theta^l$,
as in (\ref{ty}), defines a linear map 
$$\omega:\bgt^2\bbR^9\to\bgt^2\bbR^9,$$
given by
$$\bgt^2\bbR^9~\ni~ t_{ij}\quad\stackrel{\omega}{\longmapsto}\quad
\omega(t)_{kl}=\omega^{ij}_{\,\,\,\, kl}t_{ij}~\in~\bgt^2\bbR^9.$$
Here, and in the following we raise the indices by means of the
inverse $g^{ij}$ of the metric $g=g_{ij}\theta^i\theta^j$ given by
(\ref{ty}). In particular $\omega^{ij}_{\,\,\,\, kl}=g^{ip}g^{jq}\omega_{pqkl}.$

The eigenspaces of this endomorphism give the desired decomposition of
$\bgw^2\bbR^9$. We have the following proposition.
\begin{proposition}\label{p45}
The 45-dimensional vector space $\bgs^2\bbR^9$ is an $\sog(3)\times\sog(3)$ invariant
subspace in $\bgt^2\bbR^9$ corresponding to the eigenvalue $0$ of the
operator $\omega:\bgt^2\bbR^9\to\bgt^2\bbR^9$. The decomposition 
$$\bgw^2\bbR^9=V_{[2,0]}\oplus V_{[0,2]}\oplus V_{[2,4]}\oplus V_{[4,2]}$$ 
is given by:
$$
\begin{aligned}
V_{[0,2]}&=\{~\bgt^2\bbR^9\ni F_{ij}~ :~ \omega(F)_{ij}=-4 F_{ij}~\}=\soa(3)_L\\
V_{[2,0]}&=\{~\bgt^2\bbR^9\ni F_{ij}~ :~ \omega(F)_{ij}=4 F_{ij}~\}=\soa(3)_R\\
V_{[2,4]}&=\{~\bgt^2\bbR^9\ni F_{ij}~ :~ \omega(F)_{ij}=2 F_{ij}~\}\\
V_{[4,2]}&=\{~\bgt^2\bbR^9\ni F_{ij}~ :~ \omega(F)_{ij}=-2 F_{ij}~\}.
\end{aligned}
$$
The respective dimensions are 
$$\dim V_{[2,0]}=\dim V_{[0,2]}=3,\quad \dim V_{[4,2]}=\dim V_{[2,4]}=15.$$  
\end{proposition}
\begin{remark}\label{31}
Convenient bases for the 2-forms spanning $V_{[0,2]}$ and
$V_{[2,0]}$ are
$$\kappa_0^A=\tfrac12 e_{Aij}\theta^i\dz\theta^j,\quad{\rm and}\quad
\kappa_0^{A'}=\tfrac12 e_{A'ij}\theta^i\dz\theta^j.$$
Here $e_{Aij}$ and $e_{A'ij}$ are the matrix elements of the bases
$(e_A)$ and $(e_{A'})$ of $\soa(3)_L$ and
$\soa(3)_R$ as given in (\ref{dwie}). Explicitly:
\be
\begin{aligned}
-\kappa_0^1=\theta^1\dz\theta^4+\theta^2\dz\theta^5+\theta^3\dz\theta^6\\
-\kappa_0^2=\theta^1\dz\theta^7+\theta^2\dz\theta^8+\theta^3\dz\theta^9\\
-\kappa_0^3=\theta^4\dz\theta^7+\theta^5\dz\theta^8+\theta^6\dz\theta^9\\
-\kappa_0^{1'}=\theta^1\dz\theta^2+\theta^4\dz\theta^5+\theta^7\dz\theta^8\\
-\kappa_0^{2'}=\theta^1\dz\theta^3+\theta^4\dz\theta^6+\theta^7\dz\theta^9\\
-\kappa_0^{3'}=\theta^2\dz\theta^3+\theta^5\dz\theta^6+\theta^8\dz\theta^9.
\end{aligned}
\label{kappa0}\ee
Thus we have:
$$\Span_\bbR(\kappa_0^1,\kappa_0^2,\kappa_0^3)=\soa(3)_L,\quad {\rm
  and}\quad
\Span_\bbR(\kappa_0^{1'},\kappa_0^{2'},\kappa_0^{3'})=\soa(3)_R.$$
A convenient basis for the space $V_{[2,4]}$ is given
by:
\be
\begin{aligned}
&\lambda_0^1=\theta^1\dz\theta^4-\theta^3\dz\theta^6,\quad\lambda_0^2=\theta^1\dz\theta^5+\theta^2\dz\theta^4,\quad
\lambda_0^3=\theta^1\dz\theta^6+\theta^3\dz\theta^4\\
&\lambda_0^4=\theta^1\dz\theta^7-\theta^3\dz\theta^9,\quad\lambda_0^5=\theta^1\dz\theta^8+\theta^2\dz\theta^7,\quad
\lambda_0^6=\theta^1\dz\theta^9+\theta^3\dz\theta^7\\
&\lambda_0^7=\theta^2\dz\theta^5-\theta^3\dz\theta^6,\quad\lambda_0^8=\theta^2\dz\theta^6+\theta^3\dz\theta^5,
\quad\lambda_0^9=\theta^2\dz\theta^8-\theta^3\dz\theta^9\\
&\lambda_0^{10}=\theta^2\dz\theta^9+\theta^3\dz\theta^8,\quad\lambda_0^{11}=\theta^4\dz\theta^7-\theta^6\dz\theta^9,\quad
\lambda_0^{12}=\theta^4\dz\theta^8+\theta^5\dz\theta^7\\
&\lambda_0^{13}=\theta^4\dz\theta^9+\theta^6\dz\theta^7,\quad\lambda_0^{14}=\theta^5\dz\theta^8-\theta^6\dz\theta^9,\quad
\lambda_0^{15}=\theta^5\dz\theta^9+\theta^6\dz\theta^8.
\end{aligned}
\label{lambda0l}\ee
Similarly, a basis for $V_{[4,2]}$ is 
\be
\begin{aligned}
&\lambda_0^{1'}=\theta^1\dz\theta^2-\theta^7\dz\theta^8,\quad\lambda_0^{2'}=\theta^1\dz\theta^3-\theta^7\dz\theta^9,\quad
\lambda_0^{3'}=\theta^2\dz\theta^3-\theta^8\dz\theta^9\\
&\lambda_0^{4'}=\theta^1\dz\theta^5-\theta^2\dz\theta^4,\quad\lambda_0^{5'}=\theta^1\dz\theta^6-\theta^3\dz\theta^4,\quad
\lambda_0^{6'}=\theta^2\dz\theta^6-\theta^3\dz\theta^5\\
&\lambda_0^{7'}=\theta^1\dz\theta^8-\theta^2\dz\theta^7,\quad\lambda_0^{8'}=\theta^1\dz\theta^9-\theta^3\dz\theta^7,
\quad\lambda_0^{9'}=\theta^2\dz\theta^9-\theta^3\dz\theta^8\\
&\lambda_0^{10'}=\theta^4\dz\theta^5-\theta^7\dz\theta^8,\quad\lambda_0^{11'}=\theta^4\dz\theta^6-\theta^7\dz\theta^9,\quad
\lambda_0^{12'}=\theta^5\dz\theta^6-\theta^8\dz\theta^9\\
&\lambda_0^{13'}=\theta^4\dz\theta^8-\theta^5\dz\theta^7,\quad\lambda_0^{14'}=\theta^4\dz\theta^9-\theta^6\dz\theta^7,\quad
\lambda_0^{15'}=\theta^5\dz\theta^9-\theta^6\dz\theta^8.
\end{aligned}
\label{lambda0p}\ee
\end{remark}
A partial decomposition of $\bgs^2\bbR^9$ can be obtained by means of the
Casimir operator $C^{ij}_{~~kl}$ for the tensorial representation
$\otimes^2\rho$ of the irreducible representation of
$\soa(3)_L\oplus\soa(3)_R$ defined in (\ref{dwie}). To get an explicit formula
for the operator $C^{ij}_{~~kl}$ we introduce a collective index
$\mu=1,2,3,4,5,6$, so that the six vectors $(e_\mu)=(e_A,e_{A'})$ are
the basis of the Lie algebra $\soa(3)_L\oplus\soa(3)_R$. Using this basis 
one easily calculates the Killing form $k$ for
$\soa(3)_L\oplus\soa(3)_R$. We
have $$k(e_\mu,\e_\nu)=k_{\mu\nu}=-2\delta_{\mu\nu}.$$ 
The inverse of the Killing form has components
$k^{\mu\nu}=-\tfrac12\delta^{\mu\nu}.$
Then, modulo the terms proportional to the identity, the Casimir
operator $C^{ij}_{~~kl}$ reads:
$$C^{ij}_{~~kl}=k^{\mu\nu}(e_\mu^{~~i}\phantom{}_{k}e_\nu^{~~j}\phantom{}_{l}+e_\nu^{~~i}\phantom{}_{k}e_\mu^{~~j}\phantom{}_{l}).$$
Here $e_\mu^{~~i}\phantom{}_{k}$ denotes the matrix element from the $i$th raw
and $k$th column of the Lie algebra matrix $e_\mu$ given by
(\ref{dwie}). 
This defines an endomorphism 
$$C:\bgt^2\bbR^9\to\bgt^2\bbR^9$$
given by
$$\bgt^2\bbR^9~\ni~ t_{ij}\quad\stackrel{C}{\longmapsto}\quad 
C(t)_{kl}=C^{ij}_{~~~~~kl}t_{ij}~\in~\bgt^2\bbR^9.$$
We have the following proposition.     
\begin{proposition}
The Casimir operator $C$ decomposes $\bgt^2\bbR^9$ so that:
$$\bgt^2\bbR^9=V_{[0,0]}\oplus V_{[2,2]}\oplus V_{[4,4]}\oplus
W_6\oplus W_{10}\oplus W_{30}.$$
Here:
$$
\begin{aligned}
V_{[0,0]}&=\{~\bgt^2\bbR^9\ni F_{ij}~ :~ C(F)_{ij}=-4 F_{ij}~\}\\
V_{[2,2]}&=\{~\bgt^2\bbR^9\ni F_{ij}~ :~ C(F)_{ij}=-2 F_{ij}~\}\\
V_{[4,4]}&=\{~\bgt^2\bbR^9\ni F_{ij}~ :~ C(F)_{ij}=2 F_{ij}~\}\\
W_6&=\{~\bgt^2\bbR^9\ni F_{ij}~ :~ C(F)_{ij}=-3
F_{ij}~\}=V_{[2,0]}\oplus V_{[0,2]}\\
W_{30}&=\{~\bgt^2\bbR^9\ni F_{ij}~ :~ C(F)_{ij}=0~\}=V_{[2,4]}\oplus V_{[4,2]}\\
W_{10}&=\{~\bgt^2\bbR^9\ni F_{ij}~ :~ C(F)_{ij}=-
F_{ij}~\}.
\end{aligned}
$$
We further have:
$$\bgw^2\bbR^9=W_6\oplus W_{30},$$
and
$$\bgs^2\bbR^9=V_{[0,0]}\oplus V_{[2,2]}\oplus V_{[4,4]}\oplus
W_{10}.$$
The respective dimensions of the carrier spaces $W_6$, $W_{10}$ and
$W_{30}$ are: 6, 10, 30. Spaces $V_{[0,0]}$, $V_{[2,2]}$, and
$V_{[4,4]}$ have the respective dimensions 1, 9, and 25.  
\end{proposition}
The symmetric representation $W_{10}$ further decomposes onto
5-dimensional $\sog(3)\times\sog(3)$ irreducible and nonequivalent bits:
$$W_{10}=V_{[4,0]}\oplus V_{[0,4]}.$$
One can use the Casimir operator $C$ to decompose the higher rank
tensors as well. In particular, the third rank tensors,
$t_{ijk}\in \bgt^3\bbR^9$, can be decomposed using the operator
$$\tilde{C}^{ijk}_{\quad pqr}=C^{ij}_{\,\,\,\,
  pq}\delta^k_{~r}+C^{ik}_{\,\,\,\, pr}\delta^j_{~q}+C^{jk}_{\,\,\,\,
  qr}\delta^i_{~p}.$$
This defines an endomorphism 
$$\tilde{C}:\bgt^3\bbR^9\to\bgt^3\bbR^9$$
given by:
$$\bgt^3\bbR^9~\ni~ t_{ijk}\quad\stackrel{\tilde{C}}{\longmapsto}\quad 
\tilde{C}(t)_{lmn}=\tilde{C}^{ijk}_{\quad
  lmn}t_{ijk}~\in~\bgt^3\bbR^9.$$
Applying it to $\bgw^3\bbR^9$ we get:
\begin{proposition}\label{qwe}
The eigendecomposition of $\bgw^3\bbR^9$ by the operator $\tilde{C}$
is given by:
$$\bgw^3\bbR^9=Z_6\oplus Z_9\oplus Z_{30}\oplus Z_{39},$$
where
$$
\begin{aligned}
Z_6&=\{~\bgw^3\bbR^9\ni H_{ijk}~ :~ \tilde{C}(H)_{ijk}=-5
H_{ijk}~\}=V_{[2,0]}\oplus V_{[0,2]}\\
Z_9&=\{~\bgw^3\bbR^9\ni H_{ijk}~ :~ \tilde{C}(H)_{ijk}=-4
H_{ijk}~\}=V_{[2,2]}\\
Z_{30}&=\{~\bgw^3\bbR^9\ni H_{ijk}~ :~ \tilde{C}(H)_{ijk}=-2
H_{ijk}~\}=V_{[2,4]}\oplus V_{[4,2]}\\
Z_{39}&=\{~\bgw^3\bbR^9\ni H_{ijk}~ :~
\tilde{C}(H)_{ijk}=0~\}=V_{[4,4]}\oplus V_{[0,6]}\oplus V_{[6,0]}.
\end{aligned}
$$
\end{proposition}
A more refined decomposition of $\bgw^3\bbR^9$ is obtained by using the
structural 4-form $\omega$. It produces an endomorphism 
$$\tilde{\omega}:\bgw^3\bbR^9\to\bgw^3\bbR^9$$
given by:
$$\bgw^3\bbR^9~\ni~ t_{ijk}\quad\stackrel{\tilde{\omega}}{\longmapsto}\quad 
\tilde{\omega}(t)_{ijk}=3\,\omega^{lm}_{\,\,\,\,\,\,\,[ij}t_{k]lm}~\in~\bgw^3\bbR^9.$$
We have the following proposition.
\begin{proposition}\label{qwe1}
The eigendecomposition of $\bgw^3\bbR^9$ by the operator $\tilde{\omega}$
is given by:
$$\bgw^3\bbR^9=V_{[6,0]}\oplus V_{[0,6]}\oplus Z_{18}\oplus
Z_{18'}\oplus Z_{34},$$
where
$$
\begin{aligned}
V_{[0,6]}&=\{~\bgw^3\bbR^9\ni H_{ijk}~ :~ \tilde{\omega}(H)_{ijk}=-6
H_{ijk}~\}\\
V_{[6,0]}&=\{~\bgw^3\bbR^9\ni H_{ijk}~ :~ \tilde{\omega}(H)_{ijk}=6
H_{ijk}~\}\\
Z_{18}&=\{~\bgw^3\bbR^9\ni H_{ijk}~ :~ \tilde{\omega}(H)_{ijk}=4
H_{ijk}~\}=V_{[2,4]}\oplus V_{[0,2]}\\
Z_{18'}&=\{~\bgw^3\bbR^9\ni H_{ijk}~ :~
\tilde{\omega}(H)_{ijk}=-4H_{ijk}~\}=V_{[4,2]}\oplus V_{[2,0]}\\
Z_{34}&=\{~\bgw^3\bbR^9\ni H_{ijk}~ :~ \tilde{\omega}(H)_{ijk}=0~\}=V_{[2,2]}\oplus V_{[4,4]}.
\end{aligned}
$$
\end{proposition}  
Using Propositions \ref{qwe} and \ref{qwe1}, we identify all the
irreducible components of the $\sog(3)\times\sog(3)$ decomposition of
$\bgw^3\bbR^9$. For example:
$V_{[2,0]}=Z_6\cap Z_{18'}$, $V_{[4,4]}=Z_{39}\cap Z_{34}$, etc.
\section{Irreducible $\sog(3)\times\sog(3)$ geometry in dimension nine} \label{secirrgeom}
We are now prepared to define the basic object of our studies in this
article.
\begin{definition}
The irreducible $\sog(3)\times\sog(3)$ geometry in dimension nine
$(M^9,g,$ $\ten)$ is a
9-dimensional manifold $M^9$, equipped with totally
symmetric tensor fields 
$(g,\ten)$ of the respective ranks $\binom{0}{2}$ and $\binom{0}{3}$,
which at each point $x\in M^9$, reduce the structure group
$\glg(9,\bbR)$ of the tangent space ${\rm T}_xM$ to the irreducible
$(\sog(3)\times\sog(3))\subset\sog(9)\subset\glg(9,\bbR)$.\\
Alternatively, the irreducible $\sog(3)\times\sog(3)$ geometry in
dimension nine is a
9-dimensional manifold $M^9$, equipped with a differential  4-form $\omega$ which, at
each point  $x\in M^9$, 
reduces the structure group
$\glg(9,\bbR)$ of the tangent space ${\rm T}_xM$ to the irreducible
$(\sog(3)\times\sog(3))\subset\sog(9)\subset\glg(9,\bbR).$ 
\end{definition}
\begin{definition}
Given an irreducible $\sog(3)\times\sog(3)$
geometries in dimension nine $(M^9, g,\ten)$ a diffeomorphism $\phi:M^9\to
M^9$  such
that $\phi^*g=g$ and $\phi^*\ten=\ten$ is called a symmetry of
$(M^9,g,\ten)$. An infinitesimal symmetry of $(M^9,g,\ten)$ is a vector
field $X$ on $M^9$ such that ${\mathcal L}_Xg=0$ and ${\mathcal
  L}_X\ten=0$.  
\end{definition}
Symmetries of $(M^9,g,\ten)$ form a Lie group of symmetries, and
infinitesimal symmetries form a Lie algebra of symmetries. 
\subsection{$\soa(3)_L\oplus\soa(3)_R$ connection}\label{so3c} We want to analyse the properties of the irreducible
$\sog(3)\times\sog(3)$ geometries in dimension 9 by means of an 
$\soa(3)_L\oplus\soa(3)_R$-valued \emph{connection}. Since
$\soa(3)_L\oplus\soa(3)_R$ seats naturally in $\soa(9)$ such connection is
automatically \emph{metric}. It also preserves $\ten$ and
$\omega$. 

For the purpose of this paper it is convenient to think about a connection
as a Lie-algebra-valued 1-form $\Gamma$ on $M^9$. Thus, the 1-form
$\Gamma$ of the connection we are going to define for geometries
$(M^9,g,\ten,\omega)$, has values in 
$\mathfrak{g}=\soa(3)_L\oplus\soa(3)_R\subset\soa(9)$, i.e. in the Lie
algebra defined by (\ref{lieg})-(\ref{dwie}).

For further use we need the following notion:
\begin{definition}
Given an irreducible $\sog(3)\times\sog(3)$ geometry
$(M^9,g,\ten,\omega)$, a coframe
$\theta=(\theta^1,\theta^2,\theta^3,\theta^4,\theta^5,\theta^6,\theta^7,\theta^8,\theta^9)$ on $M^9$ is called adapted to it, iff
the structural tensors $g,\ten$ and $\omega$ assume
the form (\ref{ty}) in it.
\end{definition}   

Since the manifold $(M^9,g,\ten,\omega)$ is equipped with a Riemannian
metric $g$ it carries the Levi-Civita connection $\lc$ of $g$. This
can be split onto
\be\lc=\Gamma+\hbox{`the rest'}.\label{split}\ee
The only requirement that $\Gamma$ has values in $\mathfrak{g}$ is to
weak to make the above split unique. In order to achieve the uniqueness one has
to impose some (e.g. algebraic) restrictions on `the rest'. The
strongest of such restrictions is that the ${\rm `rest'}\equiv 0$. In
the next section we will provide another much weaker condition that makes
the split (\ref{split}) unique.  Here we do some preparatory steps to this.

Given the geometry $(M^9,g,\ten,\omega)$ we use a coframe $\theta$
adapted to it and write down the structure equations. This have the
form:
\be
\begin{aligned}
\der\theta^i+\Gamma^i_{~j}\dz\theta^j=T^i\\
\der\Gamma^i_{~j}+\Gamma^i_{~k}\dz\Gamma^k_{~j}=K^i_{~j}.
\end{aligned}\label{seq}
\ee
Here the matrices $\Gamma=(\Gamma^i_{~j})$ have values in the Lie
algebra $\mathfrak{g}=\soa(3)_L\oplus\soa(3)_R\subset\soa(9)$ and
therefore can be written as:
\be
\Gamma^i_{~j}=\gamma^A{e_A}^i_{~j}+\gamma^{A'}{e_{A'}}^i_{~j},\label{sex}\ee
where $(\gamma^A,\gamma^{A'})$ are 1-forms on $M^9$, and the matrices 
$e_A=({e_A}^i_{~j})$ and $e_{A'}=({e_{A'}}^i_{~j})$ are given by
(\ref{dwie}). 

The vector-valued 2-forms $$T^i=\tfrac12T^i_{~jk}\theta^j\dz\theta^k$$
represent the
`torsion' of connection $\Gamma$. The `a priori' $\soa(9)$-valued 2-forms
$$K^i_{~j}=\tfrac12K^i_{~jkl}\theta^k\dz\theta^l,$$ 
are actually $\mathfrak{g}$-valued. Hence they can also be written as
$$K^i_{~j}=\kappa^A{e_A}^i_{~j}+\kappa^{A'}{e_{A'}}^i_{~j},$$
where $$\kappa^A=\tfrac12\kappa^A_{~ij}\theta^i\dz\theta^j\quad{\rm
  and}\quad \kappa^{A'}=\tfrac12\kappa^{A'}_{~ij}\theta^i\dz\theta^j$$
are 2-forms on $M^9$. They describe the `curvature' of the connection
$\Gamma$. 

We want that the first of the structural equations (\ref{seq}), which
defines the torsion $T$ of the $\soa(3)_L\oplus\soa(3)_R$ connection
$\Gamma$, be nothing else but a reinterpretation of the `zero'-torsion equation 
\be
\der\theta^i+\lc^i_{~j}\dz\theta^j=0,\label{1stlc}\ee
for the Levi-Civita connection $\lc$. For this we need that 
$$\lc_{ijk}=\Gamma_{ijk}+\tfrac12(T_{ijk}-T_{jik}-T_{kij}),$$
or, what is the same, 
\be
\lc^i_{~j}=\Gamma^i_{~j}+\tfrac12
T^i_{~j}-\tfrac12(T_j\phantom{}^i_{~k}+T_k\phantom{}^i_{~j})\theta^k.\label{spc}\ee
Indeed, inserting the above relation into (\ref{1stlc}), because of the
symmetry of the last two terms in indices $\{jk\}$, we get
precisely the first of the structure equations (\ref{seq}).

The structural equations (\ref{seq}) when written explicitly in terms
of $(\theta^i,\gamma^A,\gamma^{A'})$ read:
\be
\begin{aligned}
\der\theta^1=\gamma^1\dz\theta^4+\gamma^2\dz\theta^7+\gamma^{1'}\dz\theta^2+\gamma^{2'}\dz\theta^3+T^1\\
\der\theta^2=\gamma^1\dz\theta^5+\gamma^2\dz\theta^8-\gamma^{1'}\dz\theta^1+\gamma^{3'}\dz\theta^3+T^2\\
\der\theta^3=\gamma^1\dz\theta^6+\gamma^2\dz\theta^9-\gamma^{2'}\dz\theta^1-\gamma^{3'}\dz\theta^2+T^3\\
\der\theta^4=-\gamma^1\dz\theta^1+\gamma^3\dz\theta^7+\gamma^{1'}\dz\theta^5+\gamma^{2'}\dz\theta^6+T^4\\
\der\theta^5=-\gamma^1\dz\theta^2+\gamma^3\dz\theta^8-\gamma^{1'}\dz\theta^4+\gamma^{3'}\dz\theta^6+T^5\\
\der\theta^6=-\gamma^1\dz\theta^3+\gamma^3\dz\theta^9-\gamma^{2'}\dz\theta^4-\gamma^{3'}\dz\theta^5+T^6\\
\der\theta^7=-\gamma^2\dz\theta^1-\gamma^3\dz\theta^4+\gamma^{1'}\dz\theta^8+\gamma^{2'}\dz\theta^9+T^7\\
\der\theta^8=-\gamma^2\dz\theta^2-\gamma^3\dz\theta^5-\gamma^{1'}\dz\theta^7+\gamma^{3'}\dz\theta^9+T^8\\
\der\theta^9=-\gamma^2\dz\theta^3-\gamma^3\dz\theta^6-\gamma^{2'}\dz\theta^7-\gamma^{3'}\dz\theta^8+T^9
\end{aligned}
\label{1str}\ee
\be
\begin{aligned}
\der\gamma^1=-\gamma^2\dz\gamma^3+\kappa^1\\
\der\gamma^2=-\gamma^3\dz\gamma^1+\kappa^2\\
\der\gamma^3=-\gamma^1\dz\gamma^2+\kappa^3\\
\der\gamma^{1'}=-\gamma^{2'}\dz\gamma^{3'}+\kappa^{1'}\\
\der\gamma^{2'}=-\gamma^{3'}\dz\gamma^{1'}+\kappa^{2'}\\
\der\gamma^{3'}=-\gamma^{1'}\dz\gamma^{2'}+\kappa^{3'}.
\end{aligned}
\label{2str}\ee
The equations (\ref{1str})-(\ref{2str}), together with their
integrability conditions implied by $\der^2\equiv 0$, encode all the geometric
information about the most general irreducible $\sog(3)\times\sog(3)$
geometry in dimension nine. They can be viewed in two ways:

\subsection{$\soa(6)$ Cartan connection} 
The standard point of view is that the equations are written just on
$M^9$. This point of view was assumed when we have introduced
(\ref{1str})-(\ref{2str}) above.

The less standard point of view is in the spirit of E. Cartan: One
considers equations (\ref{1str})-(\ref{2str}) as written on the
principal fiber bundle $$\sog(3)\times\sog(3)\to P\to M^9,$$
with the structure group $G$. This is the Cartan bundle for the
geometry $(M^9,g,\ten,\omega)$. In this point of view the (9+3+3)=15
one-forms $(\theta^i,\gamma^A,\gamma^{A'})$ are considered to live on
$P$, rather than on $M^9$. They are linearly independent at each point
of $P$ defining a prefered coframe there.   

The system may be ultimately interpreted as a system for the curvature
of a $\soa(6)$-valued Cartan connection on $P$. This connection is
defined in terms of the prefered coframe
$(\theta^i,\gamma^A,\gamma^{A'})$ on $P$ as follows. We define a
$6\times 6$ real \emph{antisymmetric} matrix 
$$\Gamma_{\rm Cartan}=\bma 0&-\gamma^1&-\gamma^2&|&\theta^1&\theta^2&\theta^3\\
                      \gamma^1&0&-\gamma^3&|&\theta^4&\theta^5&\theta^6\\
                      \gamma^2&\gamma^3&0&|&\theta^7&\theta^8&\theta^9\\-&-&-&-&-&-&-\\
                      -\theta^1&-\theta^4&-\theta^7&|&0&-\gamma^{1'}&-\gamma^{2'}\\
                       -\theta^2&-\theta^5&-\theta^8&|&\gamma^{1'}&0&-\gamma^{3'}\\
                       -\theta^3&-\theta^6&-\theta^9&|&\gamma^{2'}&\gamma^{3'}&0\ema$$
of 1-forms, and a $9\times 9$ matrix of 2-forms $K_0$
given by
$$K_0=\kappa_0^Ae_A+\kappa_0^{A'}e_{A'}.$$
The forms $(\kappa_0^A,\kappa_0^{A'})$ are the respective basis of
$\soa(3)_R$ and $\soa(3)_L$ as defined in Remark \ref{31}.
The matrix $\Gamma_{\rm Cartan}$ of 1-forms on $P$, being antisymmetric,
has values in the Lie algebra $\soa(6)$, $\Gamma_{\rm Cartan}\in
\soa(6)\otimes\bgw^1(P)$. It defines an $\soa(6)$-valued \emph{Cartan
connection} on $P$. Due to the equations (\ref{1str})-(\ref{2str}) its
\emph{curvature},
$$\tilde{R}=\der\Gamma_{\rm Cartan}+\Gamma_{\rm Cartan}\dz\Gamma_{\rm Cartan},$$
 has the form 
$$\tilde{R}=\bma 0&-R^1&-R^2&|&T^1&T^2&T^3\\
                      R^1&0&-R^3&|&T^4&T^5&T^6\\
                      R^2&R^3&0&|&T^7&T^8&T^9\\-&-&-&-&-&-&-\\
                      -T^1&-T^4&-T^7&|&0&-R^{1'}&-R^{2'}\\
                       -T^2&-T^5&-T^8&|&R^{1'}&0&-R^{3'}\\
                       -T^3&-T^6&-T^9&|&R^{2'}&R^{3'}&0\ema,$$
where
$$R^A=\kappa^A-\kappa_0^A,\quad\quad
R^{A'}=\kappa^{A'}-\kappa_0^{A'},\quad\quad A,A'=1,2,3.$$
Thus the curvature of the $\soa(6)$-Cartan connection keeps track of
both the curvature $K$ and the torsion $T$ of the
$\soa(3)_L\oplus\soa(3)_R$ connection $\Gamma$. In particular the
connection $\Gamma_{\rm Cartan}$ is \emph{flat} iff
$$T\equiv0,\quad\&\quad R\equiv 0,$$
i.e. iff the connection $\Gamma$
has \emph{vanishing} torsion, $T\equiv 0,$ and has \emph{constant
positive curvature},
$K=K_0.$
\subsection{No torsion}\label{not}
It is very easy to find all 9-dimensional irreducible
$\sog(3)\times\sog(3)$ geometries with vanishing torsion. It follows
that the system (\ref{seq}), or equivalently
(\ref{1str})-(\ref{2str}), with 
$T^i\equiv 0,$ $ i=1,2,\dots,9,$
is so rigid on $P$ that it admits only a 1-parameter family of
solutions. More specifically, the
first Bianchi identities, $\der(\der\theta^i)\equiv 0$,
$i=1,2,\dots,9$, applied to the equations (\ref{1str}), with $T^i\equiv
0$, very quickly show that the curvatures $\kappa^A$ and
$\kappa^{A'}$ must be of the form
$$\kappa^A=s \kappa_0^A,\quad\quad{\rm and}\quad\quad\kappa^{A'}=s
\kappa_0^{A'},$$
where $s$ is a real function on $P$. 
Then, the second Bianchi identities, $$\der(\der\gamma^A)\equiv 0\equiv
\der(\der\gamma^{A'}),$$ applied to (\ref{2str}) with the $\kappa's$ as
above, show that
$\der s\equiv 0,$
i.e. that the function $s$ is constant on $P$.
This proves the proposition.
\begin{proposition}
All irreducible $\sog(3)\times\sog(3)$ geometries $(M^9, g, \ten, \omega)$ 
with vanishing torsion are locally isometric to one
of the symmetric spaces
$$M^9={\mathcal G}/(\sog(3)\times\sog(3)),$$
where
$${\mathcal G} = \sog(6),\quad\quad\sog(3,3),\quad\quad{\rm or}\quad\quad(\sog(3)\times\sog(3))\rtimes_\rho\bbR^9.$$
The Riemannian metric $g$, the tensor $\ten$, and the 4-form $\omega$ 
defining the $\sog(3)\times\sog(3)$ structure are defined in
terms of the left invariant 1-forms
$(\theta^1,\theta^2,\dots,\theta^9)$, which on $P={\mathcal G}$ satisfy equations
(\ref{1str})-(\ref{2str}) and $T^i\equiv 0$. These forms, via
(\ref{ty}), define objects $g,\ten$ and $\omega$ on $P$, which descend to a well defined Riemannian metric $g$, the
symmetric tensor $\ten$ and the 4-form $\omega$ on $M^9={\mathcal
  G}/(\sog(3)\times\sog(3))$. The Levi-Civita connection of the 
metric $g$ has Einstein Ricci tensor on $M^9$, $$\rilc(g)=4sg,$$
and has holonomy reduced to $\sog(3)\times\sog(3)$. The metric $g$  is flat
if and only if $s=0$. Otherwise it is not conformally flat.
The Cartan $\soa(6)$ connection for these structures has constant
curvature, 
$$\tilde{R}=(s-1)\bma 0&-\kappa_0^1&-\kappa_0^2&|&0&0&0\\
                      \kappa_0^1&0&-\kappa_0^3&|&0&0&0\\
                      \kappa_0^2&\kappa_0^3&0&|&0&0&0\\-&-&-&-&-&-&-\\
                      0&0&0&|&0&-\kappa_0^{1'}&-\kappa_0^{2'}\\
                       0&0&0&|&\kappa_0^{1'}&0&-\kappa_0^{3'}\\
                       0&0&0&|&\kappa_0^{2'}&\kappa_0^{3'}&0\ema,$$
and is flat iff $s=1$. The symmetry group of these structures is
${\mathcal G} = \sog(6)$ for $s>0$, $\sog(3,3)$ for $s<0$ and 
$(\sog(3)\times\sog(3))\rtimes_\rho\bbR^9$ for $s=0$.

\end{proposition} 
\begin{remark}  \label{Pontriagyn}The space ${\bf SO}(6)/({\bf SO}(3)\times {\bf SO}(3))$ appearing in this proposition is just   the Grassmannian  $Gr(3,6)$ of oriented  $3$-planes in $6$-space  and  the  $4$-form  $\omega$  coincides  (up to a multiple)  with the first Pontrjagin class of the canonical $3$-plane bundle over $Gr(3,6)$ \cite{br,mic} and the $5$-form $*\omega$ is its dual. Indeed, $\omega$  is  induced by the first Pontrjagin class of the canonical $3$-plane bundle over  the Grassmannian $Gr(3,7)$. In his PhD thesis  C. Michael \cite{mic} showed that  the  $*\omega$ calibrates the special Lagrangian Grassmannian ${\bf SU}(3)/{\subset}{\bf SO}(3)\subset Gr(3,6)$ and its congruent submanifolds (and nothing else). Moreover,  he classified also  the $8$-dimensional submanifolds of $Gr(3,7)$ that are calibrated by the dual of the first Pontrjagin class  of the canonical $3$-plane bundle  (\cite{GMM}).\end{remark}
\subsection{Spin connections}

Denote by ${\mathcal C}_9$ the real Clifford algebra of the positive definite quadratic form. ${\mathcal C}_9$ is generated by the vectors of $\R^9$ and the relation
$$
v \cdot w + w \cdot v = 2 < v, w>, \quad v,w \in \R^9,
$$
holds. The spin representation of the group $\sping(9)$ is a faithful real representation in the $16$-dimensional space $\Delta_9$ of real spinors and it  is the unique irreducible  representation of the group $\sping(9)$ in dimension $16$.  With respect to this representation the orthonormal vectors $({\bf e}_1, \dots, {\bf e}_9)$ may be represented by the matrices
$$
\begin{array}{l}
{\bf e}_1 =  \sum_{k = 0}^{15} M_{16 - k, k + 1}, \quad {\bf e}_2 =  i \, \sum_{k = 0}^{15}  (-1)^k M_{16 - k, k + 1}, \\[4 pt]
 {\bf e}_3 = \sum_{k =0}^7 ( M_{15 - 2k, 2k + 1} - M_{16 - 2k, 2k + 2}), \\[4 pt]
  {\bf e}_4 =  i \, \sum_{k =0}^7 (-1)^k  (M_{15 - 2k, 2k + 1} +  M_{16 - 2k, 2k + 2}),\\[4 pt]
  {\bf e}_5 = \sum_{k =0}^3(M_{13 - 4k, 4k + 1} + M_{14 - 4 k, 4k + 2} - M_{15 - 4k, 4k + 3} - M_{16 - 4k, 4k + 4}),\\[4 pt]
  {\bf e}_6 = i \,  \sum_{k =0}^3 (-1)^ k (M_{13 - 4k, 4k + 1} + M_{14 - 4 k, 4k + 2} + M_{15 - 4k, 4k + 3} + M_{16 - 4k, 4k + 4}),\\[4 pt]
  {\bf e}_7 = \sum_{k =0}^3 (M_{9 + k, k + 1} - M_{13 + k, k + 5} + M_{1 + k, k + 9} - M_{5 + k, k + 13}),\\[4 pt]
  {\bf e}_8 =  i \,  \sum_{k =0}^7 (M_{9 + k, k + 1} - M_{1 + k, k + 9}),\\[4 pt]
  {\bf e}_9 = \sum_{k  = 0}^7 (M_{k+1,k+1} - M_{k + 9, k + 9}),
   \end{array}
$$
where by $M_{i,j}$ we denote the $16 \times 16$-matrix  having value
$1$ at its entry  $(i,j)$ and value $0$ in all the remaining
entries. In particular we have
$${\bf e}_i^2=1,\quad\quad {\bf e}_i\cdot {\bf e}_j+{\bf e}_j\cdot{\bf
  e}_i=0,\quad\quad\forall i,j=1,2,\dots,9.$$

The double covering homomorphism $\sping(9) \longrightarrow \sog(9)$ induces the isomorphism of Lie algebras $\spina(9) \longrightarrow \soa(9)$. By means of this isomorphism the basis of the Lie algebra $\spina(3)_L \oplus \spina(3)_L$ corresponding to the basis $(e_1, e_2, e_3, e'_1, e'_2, e'_3)$ of  $\soa(3)_L \oplus\soa(3)_R$ is
$$
\begin{array}{l}
{\bf E}_1 = -\frac{1}{2} ({\bf e}_1 \cdot {\bf e}_4 + {\bf e}_2\cdot {\bf e}_5 + {\bf e}_3\cdot {\bf e}_6),\\[4 pt]
{\bf E}_2 =- \frac{1}{2} ({\bf e}_1 \cdot {\bf e}_7 + {\bf e}_2
\cdot{\bf e}_8 + {\bf e}_3 \cdot {\bf e}_9),\\[4 pt]
{\bf E}_3 =- \frac{1}{2} ({\bf e}_4 \cdot {\bf e}_7 + {\bf e}_5\cdot {\bf e}_8 + {\bf e}_6\cdot {\bf e}_9),\\[4 pt]
{\bf E}'_1 =- \frac{1}{2} ({\bf e}_1  \cdot {\bf e}_2 + {\bf e}_4\cdot {\bf e}_5 + {\bf e}_6\cdot {\bf e}_8),\\[4 pt]
{\bf E}'_2 =- \frac{1}{2} ({\bf e}_1\cdot  {\bf e}_3 + {\bf e}_4\cdot {\bf e}_6 + {\bf e}_7\cdot {\bf e}_3),\\[ 4 pt]
{\bf E}'_3 =- \frac{1}{2} ({\bf e}_2\cdot  {\bf e}_3 + {\bf e}_5\cdot {\bf e}_6 + {\bf e}_8\cdot {\bf e}_9).
\end{array}
$$
Thus, in this spinorial 16-dimensional representation, we have
$$
\begin{aligned}
\spina(3)_L \oplus \spina(3)_L = \Span ({\bf E}_1, {\bf E}_2, {\bf E}_3) \oplus  \Span ({\bf E}'_1, {\bf E}'_2, {\bf E}'_3)\\  \subset  \spina(9) = \Span (\frac{1}{2} {\bf e}_i {\bf e}_j, i < j = 1, 2, \dots, 9).
\end{aligned}
$$
Now given an $\soa(3)_L\oplus\soa(3)_R$-valued connection $\Gamma=\gamma^Ae_A+\gamma^{A'}e_{A'}$ as in
(\ref{sex}), we define a spin connection 
$$\Gamma_\spina=\gamma^A{\bf E}_A+\gamma^{A'}{\bf
  E}_{A'}\in(\spina(3)_L\oplus\spina(3)_R)\otimes\bbR^9.$$
\subsection{$\soa(3)_L$ and $\soa(3)_R$ connections}
Since every $(\soa(3)_L\oplus\soa(3)_R)$-connection $\Gamma$, as
defined in Section \ref{so3c} has values in the \emph{direct sum} of
Lie algebras $\soa(3)_L$ and $\soa(3)_R$, it naturally splits onto 
$$\Gamma=\gp+\gm, \quad\quad{\rm
  with}\quad\quad\gp\in\soa(3)_L\otimes\bbR^9,
\quad{\rm and}\quad \gm\in\soa(3)_R\otimes\bbR^9. $$
Because $\soa(3)_L$ commutes with $\soa(3)_R$ this split defines two
\emph{independent} $\soa(3)$-valued connections $\gp$ and $\gm$. 
The two independent curvatures of these connections
$$\kp^i_{~j}=\der\gp^i_{~j}+\gp^i_{~k}\dz\gp^k_{~j}=\tfrac12\Rp^i_{~jkl}\theta^k\dz\theta^l$$
and
$$\km^i_{~j}=\der\gm^i_{~j}+\gm^i_{~k}\dz\gm^k_{~j}=\tfrac12\Rm^i_{~jkl}\theta^k\dz\theta^l$$
are equal to the respective
$\soa(3)_L$ and $\soa(3)_R$ parts of the curvature of $\Gamma$:
$$\Omega^i_{~j}=\der\Gamma^i_{~j}+\Gamma^i_{~k}\dz\Gamma^k_{~j}=\kp^i_{~j}+\km^i_{~j}.$$
Moreover, since,
via the identifications $\soa(3)_L=\soa(3)_L\oplus0$ and
$\soa(3)_R=0\oplus\soa(3)_R$, both $\soa(3)_L$ and $\soa(3)_R$ are
naturally included in $\soa(9)$, we can define not only the Ricci
tensor of $\Gamma$:
$$R_{ij}=R^k_{~ikj},$$
but also the corresponding Ricci tensors of $\gp$ and $\gm$:
$$\Rp_{ij}=\Rp^k_{~ikj},\quad\quad\quad  \Rm_{ij}=\Rm^k_{~ikj}.$$ 
Thus an irreducible $\sog(3)\times\sog(3)$ geometry
$(M^9,g,\ten,\omega)$ equipped with a $(\soa(3)_L\oplus\soa(3)_R)$
connection $\Gamma$ can be Einstein in several meanings:

\begin{enumerate}

\item  with respect to the Levi-Civita
connection $\lc$,
$\rilc_{ij}=\lambda g_{ij};$

\smallskip

\item with respect to the $(\soa(3)_L\oplus\soa(3)_R)$
connection $\Gamma$,
$R_{ij}=\lambda g_{ij};$

\smallskip

\item  with respect to the $\soa(3)_L$
connection $\gp$,
$\Rp_{ij}=\lambda g_{ij};$

\smallskip

\item with respect to the $\soa(3)_R$
connection $\gm$,
$\Rm_{ij}=\lambda g_{ij}.$
\end{enumerate}

Of course the functions $\lambda$ appearing in the four above
formulae, do not need to be the same.

Calculating the Ricci curvature $R_{ij}$ for the `no-torsion' examples from
Section \ref{not}, obviously yields $\rilc_{ij}=R_{ij}=4sg_{ij}$, 
since the connections $\lc$ and $\Gamma$ coincide. But it follows that in these examples
also $\gp$ and $\gm$ connections are Einstein. Actually we have
$\Rp_{ij}=\Rm_{ij}=2 sg_{ij}$ for all the examples in Section \ref{not}.  

Similar considerations as for connections $\Gamma$, $\gp$ and $\gm$,
can be performed for the spin connection $\Gamma_\spina$. Here we
have 
$$\Gamma_\spina=\ssp+\ssm, $$
with $\gp\in\spina(3)_L\otimes\bbR^9$ and
$\ssm\in\spina(3)_R\otimes\bbR^9$. Since 
$\spina(3)_L$ commutes with $\spina(3)_R$ we again have two
independent connections $\ssp$ and $\ssm$. Since they yield
essentially the same information as $\gp$ and $\gm$ we will not
comment about them anyfurther.
\section{Nearly integrable $\sog(3)\times\sog(3)$ geometries}  \label{secnearlyint}
In the previous section we discussed general $\sog(3)\times\sog(3)$
geometries in dimension nine, and general $\soa(3)_L\oplus\soa(3)_R$
connections $\Gamma$, which were obtained from the Levi-Civita
connection $\lc$ via the split (\ref{spc}). The problem with such
connections is that in general they are \emph{not} unique. In this
section we will restrict ourselves to a subclass of irreducible $\sog(3)\times\sog(3)$ geometries in
dimension nine for which the connection $\Gamma$ apearing in the formula
(\ref{spc}) will be uniquely defined. 
This class is distinguished by the following definition.
\begin{definition}
An irreducible $\sog(3)\times\sog(3)$ geometry $(M^9,g,\ten,\omega)$  is called nearly integrable iff its structural tensor
$\ten$ is a Killing tensor with respect to the Levi-Civita connection,
i.e. iff  
\be\clc_X\ten(X,X,X)=0, \quad\quad\forall X\in \hbox{T}M.\label{kt}\ee
\end{definition}
We first write the condition (\ref{kt}) in an adapted to
$(M^9,g,\ten,\omega)$ coframe $\theta$. In such a coframe we define the
Levi-Civita connection coefficients $\lc^j_{~ki}$ to be given by 
$\clc_{X_i}\theta^j=-\lc^j_{~ki}\theta^k$, where $X_i$ are 
the vector fields $X_i$ dual
on $M^9$ to the 1-forms $\theta^i$,
$X_i\hook\theta^j=\delta^j_{~i}$. The coefficients $\lc^j_{~ki}$ 
are related to the Levi-Civita connection 1-form $\lc=(\lc^i_{~j})$ via
$\lc^i_{~j}=\lc^i_{~jk}\theta^k$. In
this setting the condition (\ref{kt}) reads:
\be \lc^m_{~(ji}~\ten_{kl)m}\equiv
0.\label{nearly} \ee 
This
motivates an introduction of the map
$$\ten':\bgw^2\bbR^9\otimes\bbR^9\mapsto \bgs^4\bbR^9$$
such that \beq
\ten'(\lc)_{ijkl}&=&12\lc^p_{~(ji}~\ten_{kl)p}\nonumber\\
&=&\lc^p_{~ji}~\ten_{pkl}+\lc^p_{~ki}~\ten_{jpl}+\lc^p_{~li}~\ten_{jkp}\nonumber\\
&+&\lc^p_{~ij}~\ten_{pkl}+\lc^p_{~kj}~\ten_{ipl}+\lc^p_{~lj}~\ten_{ikp}\label{tprim}\\
&+&\lc^p_{~ik}~\ten_{pjl}+\lc^p_{~jk}~\ten_{ipl}+\lc^p_{~lk}~\ten_{ijp}\nonumber\\
&+&\lc^p_{~il}~\ten_{pjk}+\lc^p_{~jl}~\ten_{ipk}+\lc^p_{~kl}~\ten_{ijp}.\nonumber
\eeq Comparing this with (\ref{nearly}) we have the following
proposition.
\begin{proposition}\label{lop}
An irreducible $\sog(3)\times\sog(3)$ geometry $(M^9,g,\ten,\omega)$ is nearly integrable if and only if
its Levi-Civita connection $\lc\in\ker \ten'$.
\end{proposition}
It is worthwhile to note that each of the last four rows of
(\ref{tprim}) resembles the l.h.s. of the equality
$$
X^p_{~j}\ten_{pkl}+X^p_{~k}\ten_{jpl}+X^p_{~l}\ten_{jkp}=0
$$
satisfied by every matrix
$X\in\mathfrak{g}=\soa(3)_L\oplus\soa(3)_R$. Thus, $\mathfrak{g}\otimes\bbR^9\subset\ker
\ten'.$ 
Now let us consider tensors $T^i_{~jk}$, such that
$T_{ijk}=g_{il}T^l_{~jk}$ is totally antisymmetric,
$T_{ijk}=T_{[ijk]}\in\bgw^3\bbR^9$. Via $g$ we identify the space of
the  considered tensors $T^i_{~jk}$ with $\bgw^3\bbR^9$.

Because of the antisymmetry in the last pair of indices, and due 
to the first equality in (\ref{tprim}), every such $T^i_{~jk}$ also belongs
to $\ker \ten'.$ This proves the following
Lemma.
\begin{lemma}\label{ler}
Since
$$\big(\soa(3)_L\oplus\soa(3)_R\big)\otimes\bbR^9\subset\ker \ten'\quad\quad{\rm
and}\quad\quad\bgw^3\bbR^9\subset\ker \ten'$$ then
$$\Big([\big(\soa(3)_L\oplus\soa(3)_R\big)\otimes\bbR^9]+\bgw^3\bbR^9\Big)\subset\ker
\ten'.$$
\end{lemma}
It is now crucial to calculate the dimension of $\ker\ten'$. We did it
using the symbolic algebra calculation softwares Mathematica, and
independently Maple, by solving equations (\ref{nearly}) for the most
general $\lc^i_{~jk}\in \soa(9)\otimes \bbR^9$. It follows that the 
equations impose the number 186 of independent conditions on
the $\tfrac{9\times 8}{2}\times 9=324$ free coefficients
$\lc^i_{~jk}$. Thus we have:
\begin{lemma}
$$\dim\ker\ten'=324-186=138.$$
\end{lemma}
Again with the help of the Mathematica/Maple softwares we calculated the
intersection of $(\soa(3)_L\oplus\soa(3)_R)\otimes \bbR^9$ with
$\bgw^3\bbR^9$. In this way we obtained 
\begin{lemma}
$$\Big(\big(\soa(3)_L\oplus\soa(3)_R\big)\otimes \bbR^9\Big)\cap\bgw^3\bbR^9=\{0\}.$$
\end{lemma}
Comparing the dimension of
$\Big(\soa(3)_L\oplus\soa(3)_R\Big)\otimes\bbR^9$, which is 54, with the
dimension of $\bgw^3\bbR^9$, which is 84, and $\dim\ker\ten'=138$ and
using the above Lemmas, we get the following 
\begin{proposition}
\be
\ker
\ten'=\Big(\big(\soa(3)_L\oplus\soa(3)_R\big)\otimes\bbR^9\Big)\oplus\bgw^3\bbR^9.\label{ce}\ee
\end{proposition}
This leads to the following
\begin{theorem}\label{kt2}
Every nearly integrable irreducible geometry $(M^9,g,\ten,\omega)$, defines an $\soa(3)_L\oplus\soa(3)_R$-valued
connection, whose torsion is totally antisymmetric. This connection is
unique, and defined in an adapted coframe $\theta$ via the formula
\be
\lc^i_{~jk}=\Gamma^i_{~jk}+\tfrac12 T^i_{~jk},\label{cac}\ee
where $\lc^i_{~jk}$ are the Levi-Civita connection coefficients in the
coframe $\theta$, $\Gamma=(\Gamma^i_{~j})=(\Gamma^i_{~jk}\theta^k)$ is
a 1-form on $M^9$ with values in $\mathfrak{g}=\soa(3)_L\oplus\soa(3)_R$,
and $T_{ijk}=g_{il}T^l_{~jk}$ is totally antisymmetric,
i.e. $T_{ijk}=T_{[ijk]}$. 

Conversely, every irreducible $\soa(3)_L\oplus\soa(3)_R$ geometry in dimension nine
admitting a unique $\soa(3)_L\oplus\soa(3)_R$ connection with totally skew
symmetric torsion is nearly integrable.
\end{theorem}
\begin{proof} See formula (\ref{ce}) and the Proposition \ref{lop}
%
%
%
\end{proof}
\begin{definition}
The unique $\soa(3)_L\oplus\soa(3)_R$-valued connection $\Gamma$ of a nearly
integrable $\sog(3)\times\sog(3)$ geometry  $(M^9,g,\ten,\omega)$, as
described in Theorem \ref{kt2}  is called \emph{characteristic
connection} for the geometry $(M^9,g,\ten,\omega)$. 
\end{definition}
We close this section with a proposition, which relates the torsion of
the characteristic connection of a nearly integrable structure
$(M^9,g,\ten, \omega)$, and the
exterior derivatives $\der\omega$ and $\der*\omega$.
\begin{proposition}
The derivatives $\der\omega$ and $\der*\omega$ of the structural
4-forms $\omega$ and $*\omega$ of a nearly integrable geometry
$(M^9,g,\ten, \omega)$ decompose as:
\be
\der\omega\quad\in \quad V_{[2,2]}\oplus V_{[2,4]}\oplus V_{[4,2]}\oplus
V_{[4,4]},\label{deom}\ee
and
\be
\der*\omega\quad\in \quad V_{[0,2]}\oplus V_{[2,0]}\oplus V_{[0,6]}\oplus
V_{[6,0]}\oplus V_{[2,4]}\oplus
V_{[4,2]}.\label{desom}\ee
In particular, the torsion $T\in \bgw^3\bbR^9$ of the characteristic  connection is
related to these decompositions via:
$$\der\omega\equiv 0\Longleftrightarrow
 \Big(~T\in V_{[0,2]}\oplus V_{[2,0]}\oplus V_{[0,6]}\oplus V_{[6,0]}\subset\bgw^3\bbR^9~\Big),$$
and
$$\der*\omega\equiv 0\Longleftrightarrow \Big(~T\in V_{[2,2]}\oplus V_{[4,4]}\subset\bgw^3\bbR^9~\Big).$$
\end{proposition}
\begin{proof}
It follows from the 1st structure equations (\ref{1str}) that the
derivatives $\der\omega$ and $\der*\omega$ are totally expressible in
terms of the torsion components $T_{ijk}$ of the characteristic
connection. It is also clear that the relations between $\der\omega$
and $\der*\omega$ and the torsion is algebraic, and linear in the
componets of $T$. Thus each of $\der\omega$ and $\der*\omega$ must be contained
in an 84-dimensional $\sog(3)\times\sog(3)$-invariant submodule of
the respectives modules $\bgw^5\bbR^9\simeq\bgw^4 V_{[2,2]}$ and
$\bgw^6\bbR^9\simeq\bgw^3 V_{[2,2]}$. 

Now a quick calculation using Maple/Mathematica shows that the
equation $\der\omega \equiv 0$ imposes 64 conditions on the 84
components of the torsion. Similarly, one can chceck that the equation
$\der*\omega\equiv 0$ imposes 50 conditions on the torsion. Thus
$\der\omega$ has 64 independent components, and $\der*\omega$ has 50
independent components.

Comparison of these numbers with the $\sog(3)\times\sog(3)$ decompositions
of $\bgw^4 V_{[2,2]}$ and $\bgw^3 V_{[2,2]}$ given in Proposition
\ref{dec1} quickly yields to the conclusion that $\der\omega$ and
$\der*\omega$ must be in the submodules of $\bgw^5\bbR^9$ and
$\bgw^6\bbR^9$ indicated in the proposition.    

To get the decompositions (\ref{deom})-(\ref{desom}) explicitly,
dualize the forms $\der\om$ and $\der*\om$, i.e. calculate $*\der\om$,
and $*\der*\om$, and use the respective operators defined in Section
\ref{repr}.   
\end{proof}
Note that it follows from this proposition that if the torsion $T$ of
the charactersitic connection has a component in $V_{[2,4]}$, or in
$V_{[4,2]}$, then the forms $\der\omega$ and $\der*\omega$ are both nonvanishing. 
\section{Examples of nearly integrable $\sog(3)\times\sog(3)$
  geometries}
We begin this section by considering the most general situation of a \emph{nearly integrable}
irreducible geometry  $(M^9,g,\ten,\omega)$. Thus, its
characteristic connection has a general torsion in $\bgw^3\bbR^9$. 

The group $\sog(3)\times\sog(3)$ acts on the torsion space
$\bgw^3\bbR^9$ in the following way. One of the $\sog(3)$ groups in 
$\sog(3)\times\sog(3)$ is just $\exp\big(\soa(3)_L\big)$. The other is
$\exp\big(\soa(3)_R\big)$. 
Thus we have 
$$\sog(3)\times\sog(3)=\sog(3)_L\times\sog(3)_R$$
with
$$
\sog(3)_L=\exp\big(\soa(3)_L\big), \quad \quad 
\sog(3)_R=\exp\big(\soa(3)_R\big).
$$ 
The $9\times 9$ matrices $h\in\sog(3)_L$ and $h'\in\sog(3)_R$ act
on the torsion coefficients $T_{ijk}$ via:
\be
\begin{aligned}
T_{ijk}\stackrel{h}{\mapsto} (hT)_{ijk}= h^p_{~i} h^q_{~j} h^r_{~k}T_{pqr},\\
T_{ijk}\stackrel{h'}{\mapsto} (h'T)_{ijk}=h'\phantom{}^p_{~i} h'\phantom{}^q_{~j}
h'\phantom{}^r_{~k}T_{pqr}.\end{aligned}\label{adi}\ee
There is an obvious invariant of both of these actions. It is the
square of the torsion:
\be
\| T \|^2=T_{ijk}T_{pqr}g^{ip}g^{jq}g^{kr}.\label{sqt}\ee
Thus the 84-dimensional space $\bgw^3\bbR^9$ is foliated by the
$\sog(3)\times\sog(3)$-invariant 83-dimensional spheres 
$$\bbS_T=\{T_{ijk}\in\bgw^3\bbR^9~|~T_{ijk}T_{pqr}g^{ip}g^{jq}g^{kr}=r^2~\},$$
parametrized by the real parameter $r>0$. 
The group $\sog(3)\times\sog(3)$ preserves these spheres. But, for
the dimensional reasons, its action is not transitive on them. Note that
if one restricts the torsion, forcing it to lie in 
an $\sog(3)\times\sog(3)$-invariant submodule of $\bgw^3\bbR^9$, then 
the restrictions of the spheres $\bbS_T$ to this submodule will be still
invariant with respect to both actions, but the quadrics obtained
by this restriction will have smaller
dimension than 83.

For example when the torsion $T_{ijk}$ is in the
invariant module $\soa(3)_L\subset\bgw^3\bbR^9$, the spheres $\bbS_T$
restrict to 2-dimensional spheres. In such case the 3-dimensional torsion space
$\soa(3)_L\simeq\bbR^3$ is foliated by 2-dimensional spheres
with radius $r$ and center at the origin - the zero torsion. The orbit
space of the action of the groups $\sog(3)_L$ and $\sog(3)_R$ on these
spheres will be
discussed in the next subsection.   

\subsection{Torsion in $V_{[0,2]}=\soa(3)_L$}
The aim of this section is to find all \emph{nearly integrable}
irreducible geometries  $(M^9,g,\ten,\omega)$, whose
characteristic connection $\Gamma$ has totally skew symmetric torsion
$T$ in the irreducible representation $\soa(3)_L$,
$T\in\soa(3)_L\subset\bgw^3\bbR^9$.

An assumption that 
$$T\in\soa(3)_L\subset\bgw^3\bbR^9$$ 
is equivalent to the requirement, that in a coframe $\theta^i$, adapted to
$(M^9,g,\ten,\omega)$, we have
$$
\begin{aligned}
&T^i=\tfrac12 g^{ij}T_{jkl}\theta^k\dz\theta^l,\quad\quad
T_{ijk}=T_{[ijk]},\\
&\tilde{C}(T)_{ijk}=-5T_{ijk},\quad\quad{\rm and}\quad\quad
\tilde{\omega}(T)_{ijk}=4T_{ijk}.\end{aligned}$$
The last two conditions mean that, in accordance with the results of
Section \ref{repr}, the torsion is in the intersection
$Z_6\cap Z_{18}$. These algebraic conditions for $T_{ijk}$ can be
easily solved. The result is summarized in the following proposition.
\begin{proposition}\label{prtor}
In an adapted coframe $(\theta^i)$ the $\soa(3)_L$ torsion 
of  the characteristic connection of a nearly integrable 
geometry $(M^9,g,\ten,\omega)$ reads:
\be
\begin{aligned}
&\\
&T^1=-3 t_3 \theta^2\dz\theta^3+t_2 \theta^2\dz\theta^6-t_1
\theta^2\dz\theta^9-t_2\theta^3\dz\theta^5+t_1
\theta^3\dz\theta^8-t_3\theta^5\dz\theta^6-t_3 \theta^8\dz\theta^9\\&\\
&T^2=3 t_3 \theta^1\dz\theta^3-t_2 \theta^1\dz\theta^6+t_1
\theta^1\dz\theta^9+t_2\theta^3\dz\theta^4-t_1
\theta^3\dz\theta^7+t_3\theta^4\dz\theta^6+t_3 \theta^7\dz\theta^9\\&\\
&T^3=-3 t_3 \theta^1\dz\theta^2+t_2 \theta^1\dz\theta^5-t_1
\theta^1\dz\theta^8-t_2\theta^2\dz\theta^4+t_1
\theta^2\dz\theta^7-t_3\theta^4\dz\theta^5-t_3 \theta^7\dz\theta^8\\&\\
&T^4= t_2 \theta^2\dz\theta^3-t_3 \theta^2\dz\theta^6+t_3
\theta^3\dz\theta^5+3t_2\theta^5\dz\theta^6-t_1
\theta^5\dz\theta^9+t_1\theta^6\dz\theta^8+t_2 \theta^8\dz\theta^9\\&\\
&T^5=-t_2 \theta^1\dz\theta^3+t_3 \theta^1\dz\theta^6-t_3
\theta^3\dz\theta^4-3t_2\theta^4\dz\theta^6+t_1
\theta^4\dz\theta^9-t_1\theta^6\dz\theta^7-t_2 \theta^7\dz\theta^9\\&\\
&T^6= t_2 \theta^1\dz\theta^2-t_3 \theta^1\dz\theta^5+t_3
\theta^2\dz\theta^4+3t_2\theta^4\dz\theta^5-t_1
\theta^4\dz\theta^8+t_1\theta^5\dz\theta^7+t_2 \theta^7\dz\theta^8\\&\\
&T^7= -t_1 \theta^2\dz\theta^3-t_3 \theta^2\dz\theta^9+t_3
\theta^3\dz\theta^8-t_1\theta^5\dz\theta^6+t_2
\theta^5\dz\theta^9-t_2\theta^6\dz\theta^8-3t_1 \theta^8\dz\theta^9\\&\\
&T^8=t_1 \theta^1\dz\theta^3+t_3 \theta^1\dz\theta^9-t_3
\theta^3\dz\theta^7+t_1\theta^4\dz\theta^6-t_2
\theta^4\dz\theta^9+t_2\theta^6\dz\theta^7+3t_1 \theta^7\dz\theta^9\\&\\
&T^9=- t_1 \theta^1\dz\theta^2-t_3 \theta^1\dz\theta^8+t_3
\theta^2\dz\theta^7-t_1\theta^4\dz\theta^5+t_2
\theta^4\dz\theta^8-t_2\theta^5\dz\theta^7-3t_1 \theta^7\dz\theta^8.
\end{aligned}
\label{potor}\ee 
Here $(t_1,t_2,t_3)$ are the three independent components of the
torsion $T$. 
\end{proposition} 
\begin{remark}\label{prtor1}
Rewriting the above equations in terms of the basis of 2-forms
$(\kappa_0^A,\kappa_0^{A'},$ $\lambda_0^\mu,\lambda_0^{\mu'})$, as in
Remark \ref{31}, one can see that only the primed 2-forms appear above. Explicitly:
\be
\begin{aligned}
&T^1=-t_1\lambda_0^{9'}+t_2\lambda_0^{6'}+\tfrac13 t_3(5\kappa_0^{3'}-4\lambda_0^{3'}+2\lambda_0^{12'})\\
&T^2=t_1\lambda_0^{8'}-t_2\lambda_0^{5'}+\tfrac13 t_3(-5\kappa_0^{2'}+4\lambda_0^{2'}-2\lambda_0^{11'})\\
&T^3=-t_1\lambda_0^{7'}+t_2\lambda_0^{4'}+\tfrac13 t_3(5\kappa_0^{1'}-4\lambda_0^{1'}+2\lambda_0^{10'})\\
&T^4=-t_1\lambda_0^{15'}+\tfrac13
  t_2(-5\kappa_0^{3'}-2\lambda_0^{3'}+4\lambda_0^{12'})-t_3\lambda_0^{6'}\\
&T^5=t_1\lambda_0^{14'}+\tfrac13
  t_2(5\kappa_0^{2'}+2\lambda_0^{2'}-4\lambda_0^{11'})+t_3\lambda_0^{5'}\\
&T^6=-t_1\lambda_0^{13'}+\tfrac13 t_2(-5\kappa_0^{1'}-2\lambda_0^{1'}+4\lambda_0^{10'})-t_3\lambda_0^{4'}\\
&T^7=\tfrac13
  t_1(5\kappa_0^{3'}+2\lambda_0^{3'}+2\lambda_0^{12'})+t_2\lambda_0^{15'}-t_3\lambda_0^{9'}\\
&T^8=-\tfrac13
  t_1(5\kappa_0^{2'}+2\lambda_0^{2'}+2\lambda_0^{11'})-t_2\lambda_0^{14'}+t_3\lambda_0^{8'}\\
&T^9=\tfrac13 t_1(5\kappa_0^{1'}+2\lambda_0^{1'}+2\lambda_0^{10'})+t_2\lambda_0^{13'}-t_3\lambda_0^{7'}.
\end{aligned}
\label{potorl}\ee 
\end{remark}
Once the torsion in $\soa(3)_L\simeq\bbR^3$ is totally determined and
parametrized as above by a `vector' ${\bf t}=(t_1,t_2,t_3)$, we can check what are the orbits
of the action of the groups $\sog(3)_L$ and $\sog(3)_R$ on  the torsion space $\soa(3)_L\simeq\bbR^3$. A direct
calculation, yields the following two propositions:
\begin{proposition}\label{trans0}
The action of $\sog(3)_R$ on $V_{[0,2]}=\soa(3)_L$, as defined in
(\ref{adi}) 
is trivial, i.e.
 $$(h'T)_{ijk}=T_{ijk},\quad\quad\forall h'\in\sog(3)_R,\quad{\rm
  and}\quad\forall T_{ijk}\in V_{[0,2]}=\soa(3)_L.$$
\end{proposition}
On the other hand the action of $\sog(3)_L$ turns out to be as
transitive as it is only possible (remember that $\sog(3)_L$ can not
join torsions on 2-spheres $\bbS_T$ with different radii):
\begin{proposition}\label{trans}
The group $\sog(3)_L$ acts transitively on each of the 2-spheres
$\bbS_T\subset \soa(3)_L$. The orbit space of the action of
$\sog(3)_L$ on $\soa(3)\simeq\bbR^3$ is $\bbR_+\cup\{0\}$, and is 
parametrized by the radius $r$ of these spheres. Thus the orbit
structure of this action is represented by 
$$\soa(3)_L=\bbS^2\times \bbR_+\cup\{0\}.$$ 
\end{proposition}
\begin{proof}
The proof of both propositions above consists in a pure
calculation. Here we comment only on a (useful) formula for the transformation of the torsions
under the action of $\sog(3)_L$. 

Using the usual notation for the
standard scalar product of vectors ${\bf v}$ and ${\bf w}$ in
$\bbR^3$, $<{\bf v},{\bf w}>={\bf v}\cdot{\bf w}$, we anounce that 
the torsion coefficients ${\bf
  t}'=(t_1',t_2',t_3')$ transformed by $\sog(3)_L$ read:
$$
t_1'={\bf t} \cdot {\bf n}_1,\quad\quad t_2'={\bf t}\cdot {\bf n}_2,\quad\quad t_3'={\bf t}\cdot {\bf n}_3,
$$
where the vectors ${\bf n}_\mu$, $\mu=1,2,3$
are three vectors in $\bbR^3$ given by 
$$
\begin{aligned}
&{\bf n}_1=\bma\cos a_2\cos a_3\\\cos a_3\sin a_1\sin a_2+\cos a_1\sin
  a_3\\-\cos a_1\cos a_3\sin a_2+\sin a_1\sin a_3\ema,\\&\\
&{\bf n}_2=\bma-\cos a_2\sin a_3\\-\sin a_1\sin a_2\sin a_3+\cos a_1\cos
  a_3\\\cos a_1\sin a_2\sin a_3+\cos a_3\sin a_1\ema,\\&\\&{\bf n}_3=\bma\sin a_2\\ -\cos a_2\sin a_1\\\cos a_1\cos a_2\ema.
\end{aligned}
$$
They are related to a general element $h$ of the transformation group $\sog(3)_L$ via:
$$h~=~\exp(a_1e_1)\cdot\exp(a_2e_2)\cdot \exp(a_3e_3)\,\in\,\sog(3)_L,$$
where $(e_1,e_2,e_3)$ are the Lie algebra $\soa(3)_L$ generators given
by formulae (\ref{dwie}). Note that the three vectors $({\bf n}_1,{\bf
  n}_2,{\bf n}_3)$ are \emph{orthonormal}, ${\bf n}_\mu\cdot{\bf
  n}_\nu=\delta_{\mu\nu}$. 
Note also that when the group element $h$ passes 
through all the elements in $\sog(3)_L$ the three orthonormal vectors $({\bf n}_1,{\bf
  n}_2,{\bf n}_3)$ became every possible orthonormal frame
attached at the origin of $\bbR^3$. This means that given a torsion
vector ${\bf t}=(t_1,t_2,t_3)\in\soa(3)_L\simeq\bbR^3$ we can always find an element $h$ in
the group $\sog(3)_L$ which alligns the first vector ${\bf n}_1$ of the frame $({\bf n}_1,{\bf
  n}_2,{\bf n}_3)$ with ${\bf t}$. This makes 
$$t'_1=\sqrt{t_1^2+t_2^2+t_3^2},\quad\quad t'_2=0,\quad\quad t'_3=0.$$
This shows that every torsion vector ${\bf t}=(t_1,t_2,t_3)\in\soa(3)_L$ may be
  transformed to the vector $(||{\bf t}||,0,0)$. This, in particular,
  proves the transitivity of the
  $\sog(3)_L$ action on spheres with a given radius
  $T=||{\bf t}||$. 
\end{proof}

Now we analyse the differential consequences of the structure equations
(\ref{1str})-(\ref{2str}) 
with torsion $T^i$ as in (\ref{potor}). We consider the equations (\ref{1str})-(\ref{2str}) \emph{on the bundle}
$\sog(3)\times\sog(3)\to P\to M$. Thus the 15 forms
$(\theta^i,\gamma^A,\gamma^{A'})$ appearing in these equations are
considered to be linearly independent. Also the unknown torsions
$(t_1,t_2,t_3)$, as well as the curvatures, $K^i_{~jkl}$, are
considered to be functions on $P$. 

A piece of terminology is useful here: whenever we make an analysis of
a system of equations like the one given by 
(\ref{1str})-(\ref{2str}), (\ref{potor}), we will say that we analize 
an \emph{exterior differential system} - an \emph{EDS}. 

Although we have proven above that we can always gauge the
3-dimensional torsion $(t_1,t_2,t_3)$ of our EDS in such a way that
$t_2\equiv t_3\equiv 0$,
we will not use this gauge yet. This is because the use of this gauge 
would imply the
restriction of the EDS from 15-dimensional bundle $P$ to its 
13-dimensional section $P^{13}$. Since the analysis of the system is
more convenient on $P$, rather than on $P^{13}$ (because only from
there the system nicely generalizes to torsions more general than those in
$\soa(3)_L$), we will make the
gauge $t_2\equiv t_3\equiv 0$ only, after extracting the 
information from the first \emph{Bianchi identities} of our EDS on $P$. 

The first Bianchi identities are obtained by applying the exterior
derivative on
the both sides of equations (\ref{1str}). Their consequences are summarized in the
following proposition.
\begin{proposition}\label{mist}
The first Bianchi identities imply that 
\be\label{tore}
\begin{aligned}
\der t_1&=t_2\gamma^3-t_3\gamma^2\\
\der t_2&=t_3\gamma^1-t_1\gamma^3\\
\der t_3&=t_1\gamma^2-t_2\gamma^1,
\end{aligned}
\ee
and that the curvatures $(\kappa^A,\kappa^{A'})$, as defined in
(\ref{2str}), read:
\be\label{curve}
\begin{aligned}
&\kappa^1~=~k~\kappa_0^1~+~t_1t_2~\kappa_0^2~+~t_1t_3~\kappa_0^3\\
&\kappa^2~=~t_1t_2~\kappa_0^1~+~(k-t_1^2+t_2^2)~\kappa_0^2~+~t_2t_3~\kappa_0^3\\
&\kappa^3~=~t_1t_3~\kappa_0^1~+~t_2t_3~\kappa_0^2~+~(k-t_1^2+t_3^2)~\kappa_0^3\\
&\kappa^{1'}~=~(k+t_1^2+2t_2^2+2t_3^2)~\kappa_0^{1'}\\
&\kappa^{2'}~=~(k+t_1^2+2t_2^2+2t_3^2)~\kappa_0^{2'}\\
&\kappa^{3'}~=~(k+t_1^2+2t_2^2+2t_3^2)~\kappa_0^{3'},
\end{aligned}
\ee
Here $k$ is an unknown function on $P$, and the forms
$(\kappa_0^A,\kappa_0^{A'})$ are defined in (\ref{kappa0}). 
\end{proposition} 
Thus, the first Bianchi identities show that the curvature of the
characteristic connection is totally determined by the torsion
$(t_1,t_2,t_3)$ and an unknown function $k$. 
\begin{proof} (of the Proposition). To apply the first Bianchi
  identities, one needs the derivatives of the torsions $t_i$. So we
  assume the most general form for these:
\be
\der
  t_\mu=t_{\mu j}\theta^j+t_{\mu A}\gamma^A+t_{\mu A'}\gamma^{A'},\quad\quad
  \mu=1,2,3.\label{pom}\ee
Here $t_{\mu j}, t_{\mu A}, t_{\mu A'}$ are (3*9+3*3+3*3)=45 functions on $P$, which we hope to
determine by means of the first Bianchi identities
$\der^2\theta^i\equiv 0$, $i=1,2,\dots,9$. Note that if one applies
the exterior differential to the equations (\ref{1str}), the $\der$ of the right hand
sides must be zero, $\der({\rm rhs})\equiv 0$. Inserting our
definitions (\ref{pom}) in these identities, we obtain nine identities
each of which is a 3-form on $P$. Decomposing these nine 3-forms onto
the basis of 3-forms on $P$, which consists of the primitive forms 
$\theta^i\dz\theta^j\dz\theta^k$,
$\theta^i\dz\theta^j\dz\gamma^{A/A'}$,
$\theta^i\dz\gamma^{A/A'}\dz\gamma^{B/B'}$, and
$\gamma^{A/A'}\dz\gamma^{B/B'}\dz\gamma^{C/C'}$, one gets relations on
the unknown functions $t_{\mu j}, t_{\mu A}, t_{\mu A'}$ and the
curvatures $K^i_{~jkl}$. 
  
Analysing these relations step by
step we get the following:
\begin{itemize}
\item First, we consider terms at the basis forms
  $\theta^i\dz\theta^j\dz\gamma^{A/A'}$. This gives 18 conditions
  determining all the functions $t_{\mu A}$ and $t_{\mu A'}$ in terms
  of $(t_1,t_2,t_3)$. After solving these 18 conditions we get:
$$
\begin{aligned}
\der t_1&=t_2\gamma^3-t_3\gamma^2+t_{1j}\theta^j\\
\der t_2&=t_3\gamma^1-t_1\gamma^3+t_{2j}\theta^j\\
\der t_3&=t_1\gamma^2-t_2\gamma^1+t_{3j}\theta^j.
\end{aligned}
$$  
\item Second, the terms at the basis forms
  $\theta^i\dz\theta^j\dz\theta^k$ when equated to zero, can be split
  into two types of equations. The first type is obtained by
  eliminating the curvatures $K^i_{~jkl}$ from the full set. This
  yields a system of linear equations for the unknowns $t_{\mu j}$,
  whose only solution is $t_{\mu j}=0.$
After these conditions are imposed the second type of equations,
involves the curvatures $K^i_{~jkl}$ only in a linear fashion. It has a
unique solution for the curvatures, which explicitly is given by
(\ref{curve}).
\item Third, after imposing the conditions described above,
  all the other terms in $\der^2\theta^i$ are automatically zero. 
\end{itemize}  
This proves the proposition, and also shows that the conditions
(\ref{tore})-(\ref{curve}) on
the curvature and the derivatives of the torsion are
equivalent to the first Bianchi identities of the system in consideration. 
\end{proof}

Now we are in a position to impose the gauge
\be
t_2\equiv t_3\equiv 0.\label{gau}
\ee
Proposition \ref{trans} guarantees that every nearly integrable
$\sog(3)\times\sog(3)$ geometry with torsion in $\soa(3)_L$ admits an
adapted frame in which the conditions (\ref{gau}) hold. But the
asumption of the gauge (\ref{gau}) reduces the degrees of freedom
by 2, from 15 to 13. This means that we reduce the equation of our
EDS (\ref{1str})-(\ref{2str}), (\ref{potor})
from dimension 15 to dimension 13. Also the differential consequences
(\ref{tore})-(\ref{curve}) of this EDS must be reduced to dimension
13. This in particular means that the fifteen 1-forms
$(\theta^i,\gamma^A,\gamma^{A'})$ can no longer be
linearly independent. This obvious observation finds its confirmation in the integrability
conditions (\ref{tore})-(\ref{curve}). 

Indeed, assuming $t_2\equiv t_3\equiv 0$, and comparing it with the last two
integrabilty conditions (\ref{tore}) yields:
$$t_1 \gamma^3\equiv 0,\quad\quad{\rm and}\quad\quad t_1\gamma^2\equiv 0.$$
These, when confronted with the assumption that the torsion $T^i$ is
not vanishing in a neighbourhood, implies that 
\be
\gamma^2\equiv 0,\quad\quad{\rm and}\quad\quad\gamma^3\equiv 0.\label{gac}\ee
Thus the EDS (\ref{1str})-(\ref{2str}), (\ref{potor}) naturally
reduces to 13-dimensions, and has now thirteen 1-forms
$(\theta^i,\gamma^1,\gamma^{A'})$ linearly independent at each point
of the 13-dimensional manifold, which we previously called $P^{13}$. 

The relations
(\ref{gac}) have further consequences, for if we compare them with the
second and the third equation (\ref{2str}), we see that 
$$\kappa^2\equiv 0,\quad\quad{\rm and}\quad\quad\kappa^3\equiv 0.$$
If we now compare these with (\ref{gac}), and the second and the third
of integrability conditions (\ref{curve}), we get:
$$(k-t_1^2)~\kappa_0^2\equiv 0,\quad\quad{\rm and}\quad\quad
(k-t_1^2)~\kappa_0^3\equiv 0.$$
These hold iff
$$k\equiv t_1^2,$$
which we have to accept form now on. Note that this totally determines
the function $k$, which was a misterious unknown in Proposition
\ref{mist}. 

Finally, if we insert $t_2\equiv t_3\equiv 0$ in the first of the
integrabilty conditions (\ref{tore}), we get also that
$$\der t_1\equiv 0,$$
i.e. that the function $t_1$ must be \emph{constant} on the
13-dimensional reduced manifold $P^{13}$ on which our EDS lives. 

These considerations, when compared with the rest of the integrability
conditions (\ref{curve}), prove the following proposition.
\begin{proposition}
Every nearly integrable $\sog(3)\times\sog(3)$ geometry  $(M^9,g,\ten,\omega)$ with a 
nonvanishing torsion $T$ of the characteristic connection lying in
$\soa(3)_L=V_{[0,2]}$, $T\in\soa(3)_L$,  can be
described in terms of thirteen linearly independent 1-forms
$(\theta^i,\gamma^1,\gamma^{A'})$, $i=1,2,\dots, 9$, $A'=1,2,3$,
satisfying 
\be
\begin{aligned}
&\der\theta^1=\gamma^1\dz\theta^4+\gamma^{1'}\dz\theta^2+\gamma^{2'}\dz\theta^3+t~(-\theta^2\dz\theta^9+\theta^3\dz\theta^8)\\
&\der\theta^2=\gamma^1\dz\theta^5-\gamma^{1'}\dz\theta^1+\gamma^{3'}\dz\theta^3+t~(\theta^1\dz\theta^9-\theta^3\dz\theta^7)\\
&\der\theta^3=\gamma^1\dz\theta^6-\gamma^{2'}\dz\theta^1-\gamma^{3'}\dz\theta^2+t~(-\theta^1\dz\theta^8+\theta^2\dz\theta^7)\\
&\der\theta^4=-\gamma^1\dz\theta^1+\gamma^{1'}\dz\theta^5+\gamma^{2'}\dz\theta^6+t~(-\theta^5\dz\theta^9+\theta^6\dz\theta^8)\\
&\der\theta^5=-\gamma^1\dz\theta^2-\gamma^{1'}\dz\theta^4+\gamma^{3'}\dz\theta^6+t~(\theta^4\dz\theta^9-\theta^6\dz\theta^7)\\
&\der\theta^6=-\gamma^1\dz\theta^3-\gamma^{2'}\dz\theta^4-\gamma^{3'}\dz\theta^5+t~(-\theta^4\dz\theta^8+\theta^5\dz\theta^7)\\
&\der\theta^7=\gamma^{1'}\dz\theta^8+\gamma^{2'}\dz\theta^9-t~(\theta^2\dz\theta^3+\theta^5\dz\theta^6+3
  \theta^8\dz\theta^9)\\
&\der\theta^8=-\gamma^{1'}\dz\theta^7+\gamma^{3'}\dz\theta^9+t~(\theta^1\dz\theta^3+\theta^4\dz\theta^6+3
  \theta^7\dz\theta^9)\\
&\der\theta^9=-\gamma^{2'}\dz\theta^7-\gamma^{3'}\dz\theta^8-t~(\theta^1\dz\theta^2+\theta^4\dz\theta^5+3
  \theta^7\dz\theta^8)
\end{aligned}
\label{131}\ee
\be
\begin{aligned}
&\der\gamma^1=~t^2~(\theta^1\dz\theta^4+\theta^2\dz\theta^5+\theta^3\dz\theta^6)\\
&\der\gamma^{1'}=-\gamma^{2'}\dz\gamma^{3'}+2t^2~(\theta^1\dz\theta^2+\theta^4\dz\theta^5+\theta^7\dz\theta^8)\\
&\der\gamma^{2'}=-\gamma^{3'}\dz\gamma^{1'}+2t^2~(\theta^1\dz\theta^3+\theta^4\dz\theta^6+\theta^7\dz\theta^9)\\
&\der\gamma^{3'}=-\gamma^{1'}\dz\gamma^{2'}+2t^2~(\theta^2\dz\theta^3+\theta^5\dz\theta^6+\theta^8\dz\theta^9).
\end{aligned}
\label{132}\ee
Here $\der t\equiv 0$, i.e. the function $t$ is constant.
\end{proposition}  

Note, that the system (\ref{131})-(\ref{132}) involves only constant
coefficients on the right hand sides. Thus the manifold $P^{13}$ is a
Lie group $P^{13}={\mathcal G}^{13}$, with the forms
$(\theta^i,\gamma^1,\gamma^{A'})$ constituting a basis of its left invariant forms. 
A calculation of the Killing form for ${\mathcal G}^{13}$, by using the structure
constants red off from (\ref{131})-(\ref{132}), shows that
this group is semisimple, unless the torsion $t\equiv 0$. The group
${\mathcal G}^{13}$ is a transitive group of symmetries of the underlying nearly
integrable geometry $(M^9,g,\ten,\omega)$. The 9-dimensional manifold
$M^9$ is a homogeneous space $M^9={\mathcal G}^{13}/H$, where $H$ is a certain
4-dimensional subgroup of ${\mathcal G}^{13}$. 
The structural tensors $g$, $\ten$
and $\omega$ of the corresponding $\sog(3)\times\sog(3)$ structure are
obtained, via formulae (\ref{ty}), from the 1-forms $(\theta^i)$ solving
(\ref{131})-(\ref{132}). The system
(\ref{131})-(\ref{132}) guarantees that although tensors 
$g,\ten,\omega$ defined in this way live on ${\mathcal G}^{13}$, they actually descend to
tensors $g,\ten,\omega$ on the manifold
$M^9={\mathcal G}^{13}/H$, defining a homogeneus nearly integrable geometry $(M^9, g,\ten,\omega)$ with
13-dimensional group of symmetries ${\mathcal G}^{13}$ there. 

  For $t =0$ the Lie group ${\mathcal G}^{13}$ is just  a semidirect product $(\sog(3) \times \sog(2)) \ltimes \R^9$. For $t \neq 0$,  by considering the new basis of $1$-forms 
$$
\begin{array}{l}
\tilde \theta^i =  t \theta^ i, i = 1, \ldots, 6,\\
\tilde \gamma^1 = \gamma'^1+ t \theta^9,   \quad
\tilde \gamma^2 = \gamma'^2 - t \theta^8,  \quad
\tilde \gamma^3 = \gamma'^3 + t \theta^7, \quad 
\\
\tilde \theta^7 = \gamma'^3 + 2 t \theta^7, \quad
\tilde \theta^8 = \gamma'^2 - 2 t \theta^8, \quad
\tilde \theta^9 = \gamma'^1+ 2 t \theta^9,
\end{array}
$$  
one sees that for any $t \neq 0$ the Lie group ${\mathcal G}^{13}$  is  the product $\sog(3) \times K^{10}$  with structure equations
$$
\begin{array}{l}
d  \tilde \theta^1 = \gamma^1 \wedge \tilde \theta^4 +  \tilde \gamma^1 \wedge \tilde \theta^2 + \tilde \gamma^2 \wedge \tilde \theta^3,\\
d  \tilde \theta^2 = \gamma^1 \wedge \tilde \theta^5 -  \tilde \gamma^1 \wedge \tilde \theta^1+ \tilde \gamma^3 \wedge \tilde \theta^3,\\
d  \tilde \theta^3 = \gamma^1 \wedge \tilde \theta^6 -  \tilde \gamma^2 \wedge \tilde \theta^1 + \tilde \gamma^2 \wedge \tilde \theta^2,\\
d  \tilde \theta^4 = - \gamma^1 \wedge \tilde \theta^1 +  \tilde \gamma^1 \wedge \tilde \theta^5 + \tilde \gamma^2 \wedge \tilde \theta^6,\\
d  \tilde \theta^5 = - \gamma^1 \wedge \tilde \theta^2 -  \tilde \gamma^1 \wedge \tilde \theta^4 + \tilde \gamma^3 \wedge \tilde \theta^6,\\
d  \tilde \theta^6 = - \gamma^1 \wedge \tilde \theta^3 -  \tilde \gamma^2 \wedge \tilde \theta^4 - \tilde \gamma^3 \wedge \tilde \theta^5,\\
 d \gamma^1 =  \tilde \theta^{1} \wedge  \tilde \theta^4 + \tilde  \theta^2 \wedge  \tilde \theta^5 + \tilde  \theta^3 \wedge \tilde \theta^6,\\
 d\tilde \gamma^1 =  - \tilde \gamma^2 \wedge \tilde \gamma^3 + \tilde \theta^1 \wedge  \tilde \theta^2 + \tilde \theta^4 \wedge \tilde \theta^5,\\
d\tilde \gamma^2 =  - \tilde \gamma^3 \wedge \tilde \gamma^1 + \tilde \theta^1 \wedge  \tilde \theta^3 + \tilde \theta^4 \wedge \tilde \theta^6,\\
d\tilde \gamma^3 =  - \tilde \gamma^1\wedge \tilde \gamma^2 + \tilde \theta^2 \wedge \tilde  \theta^3 + \tilde \theta^5 \wedge \tilde \theta^6,\\
d  \tilde \theta^7 =   \tilde \theta^8 \wedge  \tilde \theta^9,\\
d  \tilde \theta^8 =   \tilde \theta^9 \wedge   \tilde \theta^7,\\
d  \tilde \theta^9 =   \tilde \theta^7 \wedge  \tilde \theta^8.\\
\end{array}
$$
To say what is $K^{10}$ we calculate the Killing forms. In the basis
$(\tilde\theta^1,\tilde\theta^2,\tilde\theta^3,\tilde\theta^4,\tilde\theta^5,\tilde\theta^6,\gamma^1,$
$\tilde{\gamma}^1,\tilde{\gamma}^2,\tilde{\gamma}^3,\tilde\theta^7,\tilde\theta^8,\tilde\theta^9)$
the Killing form of ${\mathcal G}^{13}$ reads:
$$Kil_{13}={\rm diag}(6,6,6,6,6,6,-6,-6,-6,-6,-2,-2,-2).$$
The Lie algebra of $K^{10}$ is spanned by
$(\tilde\theta^1,\tilde\theta^2,\tilde\theta^3,\tilde\theta^4,\tilde\theta^5,\tilde\theta^6,\gamma^1,\tilde{\gamma}^1,\tilde{\gamma}^2,\tilde{\gamma}^3)$. Its
Killing form in this basis is:
$$Kil_{10}={\rm diag}(6,6,6,6,6,6,-6,-6,-6,-6),$$
showing that $K^{10}$ is semisimple, and as such, having 
dimension 10, it must
be locally isomorphic to a noncompact real form of $\sog(5,\bbC)$. Comparison of
Killing forms for $\sog(1,4)$ and $\sog(2,3)$ shows that $K^{10}$ is
locally $\sog(2,3)$.

In both cases  ($t = 0$ and $t \neq 0$) the Lie algebra of the 
group $H=\sog(3) \times \sog(2)$  is given by the  annihilator of the 1-forms $\theta^i$ , $i = 1, 2, . . . , 9$.

After calculating the curvatures of  the various connections associated
with this geometry we get the following theorem.

\begin{theorem}\label{v02}
Every nearly integrable irreducible $\sog(3)\times\sog(3)$ geometry
$(M^9,g, $ $\ten, \omega)$ with torsion of the characteristic connection
$\Gamma$ in $V_{[0,2]}=\soa(3)_L$ is locally a homogeneous space
${\mathcal G}^{13}/H$. It
has a transitive symmetry group ${\mathcal G}^{13}$ of dimension
13. For $t=0$ the  Lie group ${\mathcal G}^{13}$ is a semirect product
$(\sog(3)  \times \sog(2))\ltimes \R^9$ and for $t \neq 0$ it is a
direct product $\sog(3)\times\sog(2,3)$.

The metric $g$ is conformally \emph{non}-flat and \emph{not} locally symmetric. The Ricci tensors of the Levi-Civita
connection $\lc$, of the characteristic connection $\Gamma$, and of the
$\soa(3)_L$ part $\gp$ of the characteristic connection have all two 
distinct eigenvalues. 

The $\soa(3)_R$ part $\gm$ of the characteristic connection is
\emph{Einstein}. 

Explicitly, in the adapted coframe $(\theta^i)$ in which the structure
equations read as in (\ref{131}) and in which the structural tensors
$g,\ten,\omega$ are given by (\ref{ty}), we have:
\begin{itemize}
\item 
The Cartan connection $\Gamma_{\rm Cartan}$ has the curvature given by:
$$\tilde{R}=\bma 0&(1+t^2)\kappa_0^1&\kappa_0^2&|&T^1&T^2&T^3\\
                      -(1+t^2)\kappa_0^1&0&\kappa_0^3&|&T^4&T^5&T^6\\
                      -\kappa_0^2&-\kappa_0^3&0&|&T^7&T^8&T^9\\-&-&-&-&-&-&-\\
                      -T^1&-T^4&-T^7&|&0&(1+2t^2)\kappa_0^{1'}&(1+2t^2)\kappa_0^{2'}\\
                       -T^2&-T^5&-T^8&|&-(1+2t^2)\kappa_0^{1'}&0&(1+2t^2)\kappa_0^{3'}\\
                       -T^3&-T^6&-T^9&|&-(1+2t^2)\kappa_0^{2'}&-(1+2t^2)\kappa_0^{3'}&0\ema,$$
where the torsions $T^i$ are given by (\ref{potorl}) with $t_1=t={\rm const},
t_2=t_3=0$. 
\item The Levi-Civita connection Ricci tensor reads:
$$\rilc={\rm
  diag}\Big(-4t^2,-4t^2,-4t^2,-4t^2,-4t^2,-4t^2,~\tfrac32t^2,~\tfrac32t^2,~\tfrac32t^2\Big),$$
and has the Ricci scalar equal to $-\tfrac{39}{2}t^2.$
\item The $\soa(3)_L$ part $\gp$ of the characteristic connection has
  the curvature
$$\kp=-t^2\kappa_0^1e_1,$$
where the matrix $e_1=(e_1\phantom{}^i_j)$ is given by
(\ref{dwie}). It has the Ricci tensor $\Rp_{ij}$ given by
$$\Rp_{ij}={\rm
  diag}\Big(-t^2,-t^2,-t^2,-t^2,-t^2,-t^2,~0,~0,~0\Big),$$
with the Ricci scalar equal to $-6t^2$. 
\item The $\soa(3)_R$ part $\gm$ of the characteristic connection has
  the curvature
$$\km=-2t^2\kappa_0^{A'}e_{A'},$$
where as before the matrices $e_{A'}=(e_{A'}\phantom{}^i_j)$ are given by
(\ref{dwie}). Its Ricci tensor is \emph{Einstein}
$$\Rm_{ij}=-4 t^2g_{ij},$$
and has Ricci scalar equal to $-36 t^2$. 
\item The characteristic connection $\Gamma=\gp+\gm$ has curvature
$$\Omega=\kp+\km=-t^2\kappa_0^1e_1-2t^2\kappa_0^{A'}e_{A'}$$
and the Ricci
  tensor 
$$R_{ij}={\rm 
diag}\Big(-5t^2,-5t^2,-5t^2,-5t^2,-5t^2,-5t^2,-4t^2,-4t^2,-4t^2\Big).$$
 \end{itemize}
\end{theorem}

\subsection{Torsion in $V_{[0,6]}$}\label{sv06}
Now we find examples of \emph{nearly integrable} geometries  $(M^9,g$, $\ten,\omega)$ in dimension nine, whose
characteristic connection $\Gamma$ has totally skew symmetric torsion
$T$ in the irreducible representation $V_{[0,6]}$,
$T\in V_{[0,6]}\subset\bgw^3\bbR^9$.

The  assumption that 
$T\in V_{[0,6]}\subset\bgw^3\bbR^9$
is equivalent to the requirement, that in a coframe $\theta^i$, adapted to
$(M^9,g,\ten,\omega)$, we have
$$
T^i=\tfrac12 g^{ij}T_{jkl}\theta^k\dz\theta^l,\quad\quad
T_{ijk}=T_{[ijk]},\quad{\rm and}\quad\tilde{\omega}(T)_{ijk}=-6T_{ijk}.
$$
Solving these algebraic conditions for $T_{ijk}$ we get the following
proposition.
\begin{proposition}\label{prtor7}
In an adapted coframe $(\theta^i)$ the $V_{[0,6]}$ torsion 
of a characteristic connection of a nearly integrable 
geometry $(M^9,g,\ten,\omega)$ reads:
\be
\begin{aligned}
&T^1=u_1(-\lambda_0^{3'}+\lambda_0^{12'})-u_2\lambda_0^{15'}-u_3\lambda_0^{3'}-u_4\lambda_0^{6'}-u_5\lambda_0^{9'}-u_6\lambda_0^{6'}-u_7\lambda_0^{9'}\\
&T^2=u_1(\lambda_0^{2'}-\lambda_0^{11'})+u_2\lambda_0^{14'}+u_3\lambda_0^{2'}+u_4\lambda_0^{5'}+u_5\lambda_0^{8'}+u_6\lambda_0^{5'}+u_7\lambda_0^{8'}\\
&T^3=u_1(-\lambda_0^{1'}+\lambda_0^{10'})-u_2\lambda_0^{13'}-u_3\lambda_0^{1'}-u_4\lambda_0^{4'}-u_5\lambda_0^{7'}-u_6\lambda_0^{4'}-u_7\lambda_0^{7'}\\
&T^4=u_1\lambda_0^{6'}-u_2\lambda_0^{9'}+u_4(-\lambda_0^{3'}+\lambda_0^{12'})+u_5\lambda_0^{15'}-u_6\lambda_0^{3'}\\
&T^5=-u_1\lambda_0^{5'}+u_2\lambda_0^{8'}+u_4(\lambda_0^{2'}-\lambda_0^{11'})-u_5\lambda_0^{14'}+u_6\lambda_0^{2'}\\
&T^6=u_1\lambda_0^{4'}-u_2\lambda_0^{7'}+u_4(-\lambda_0^{1'}+\lambda_0^{10'})+u_5\lambda_0^{13'}-u_6\lambda_0^{1'}\\
&T^7=-u_2\lambda_0^{6'}+u_3\lambda_0^{9'}+u_5(-\lambda_0^{3'}+\lambda_0^{12'})+u_6\lambda_0^{15'}-u_7\lambda_0^{3'}\\
&T^8=u_2\lambda_0^{5'}-u_3\lambda_0^{8'}+u_5(\lambda_0^{2'}-\lambda_0^{11'})-u_6\lambda_0^{14'}+u_7\lambda_0^{2'}\\
&T^9=-u_2\lambda_0^{4'}+u_3\lambda_0^{7'}+u_5(-\lambda_0^{1'}+\lambda_0^{10'})+u_6\lambda_0^{13'}-u_7\lambda_0^{1'},
\end{aligned}\label{mnb}
\ee 
where $(u_1,u_2,u_3,u_4,u_5,u_6,u_7)$ are the seven independent components of the
torsion $T$, and $(\lambda_0^{\mu'})$, $\mu'=1,2,\dots, 15$, is a basis of 2-forms in
$V_{[4,2]}$ as defined in (\ref{lambda0p}). 
\end{proposition} 
Now we have an analog of Proposition \ref{trans0}:
\begin{proposition}
The action of $\sog(3)_R$ on $V_{[0,6]}$, as defined in (\ref{adi}), 
is trivial, i.e.
 $$(h'T)_{ijk}=T_{ijk},\quad\quad\forall h'\in\sog(3)_R,\quad{\rm
  and}\quad\forall T_{ijk}\in V_{[0,6]}.$$
\end{proposition}
The `left' $\sog(3)$ acts nontrivially on $V_{[0,6]}$. It has a
4-parameter family of generic orbits in this 7-dimensional space. As
in the $V_{[0,2]}$ case, instead of restricting ourselves to the
representatives of these orbits, we will analyze the EDS
(\ref{1str})-(\ref{2str}) for the torsion in $V_{[0,6]}$, with general
torsions $(u_1,u_2,u_3,u_4,u_5,u_6,u_7)$ as in
(\ref{mnb}). Thus the EDS (\ref{1str})-(\ref{2str}), (\ref{mnb}) we
consider, lives on the Cartan bundle $\sog(3)_L\times\sog(3)_R\to P\to M$,
where the 15 forms $(\theta^i,\gamma^A,\gamma^{A'})$ are linearly
independent at each point.  

Now the $V_{[0,6]}$ analog of Proposition \ref{mist} reads:
\begin{proposition}\label{koe}
The first Bianchi identities $\der^2\theta^i\equiv 0$, for the  
EDS (\ref{1str})-(\ref{2str}), (\ref{mnb}) imply that:
\be\label{toreu}
\begin{aligned}
&\der u_1=(3u_4+2u_6)\gamma^1+u_5\gamma^2-2u_2\gamma^3\\
&\der u_2=-(2u_5+u_7)\gamma^1-(u_4+2u_6)\gamma^2+(u_1-u_3)\gamma^3\\
&\der u_3=u_6\gamma^1+(2u_5+3u_7)\gamma^2+2u_2\gamma^3\\
&\der u_4=-3u_1\gamma^1+3u_5\gamma^3\\
&\der u_5=2u_2\gamma^1-u_1\gamma^2+(2u_6-u_4)\gamma^3\\
&\der u_6=-u_3\gamma^1+2u_2\gamma^2+(u_7-2u_5)\gamma^3\\
&\der u_7=-3u_3\gamma^2-3u_6\gamma^3,
\end{aligned}
\ee
and that the curvatures $(\kappa^A,\kappa^{A'})$, as defined in
(\ref{2str}), are:
\be\label{curveu}
\begin{aligned}
&\kappa^1~=~k_1~\kappa_0^1~+~k_2~\kappa_0^2~+~k_3~\kappa_0^3\\
&\kappa^2~=~k_2~\kappa_0^1~+~k_4~\kappa_0^2~+~k_5~\kappa_0^3\\
&\kappa^3~=~k_3~\kappa_0^1~+~k_5~\kappa_0^2~+~k_6~\kappa_0^3\\
&\kappa^{1'}~=~k_7~\kappa_0^{1'}\\
&\kappa^{2'}~=~k_7~\kappa_0^{2'}\\
&\kappa^{3'}~=~k_7~\kappa_0^{3'},
\end{aligned}
\ee
where:
\be\label{curveuu}
\begin{aligned}
&k_2=2(u_1+u_3)u_2-(2u_4+3u_6)u_5-(u_4+2u_6)u_7\\
&k_3=2u_2u_4+(2u_1-u_3)u_5+u_1u_7\\
&k_4=k_1+2u_1^2-2u_3^2+2u_4^2+2u_4u_6-2u_5u_7-2u_7^2\\
&k_5=-u_3u_4+(u_1-2u_3)u_6-2u_2u_7\\
&k_6=k_1+2u_1^2+2u_1u_3+2u_4^2+4u_4u_6+2u_5u_7\\
&k_7=k_1+2u_1^2+u_2^2+u_1u_3+2u_4^2+u_5^2+3u_4u_6+u_6^2+u_5u_7.
\end{aligned}
\ee
Here $k_1$ is an unknown function, and $(\kappa_0^A,\kappa_0^{A'})$
are given by (\ref{kappa0}). 
\end{proposition}
\begin{proof}
The proof here is very similar to the proof of Proposition
\ref{mist}. So we first assume the most general form for the derivatives of the torsions $u_\mu$:
\be
\der
  u_\mu=u_{\mu j}\theta^j+u_{\mu A}\gamma^A+u_{\mu A'}\gamma^{A'},\quad\quad
  \mu=1,2,\dots 7.\label{pomu}\ee
Here $u_{\mu j}, u_{\mu A}, u_{\mu A'}$ are (7*9+7*3+7*3)=105
functions on $P$, which we will 
determine by means of the first Bianchi identities
$\der^2\theta^i\equiv 0$, $i=1,2,\dots,9$. Inserting our
definitions (\ref{pomu}) in these identities, we obtain nine identities
each of which is a 3-form on $P$. We decompose these nine 3-forms onto
the basis of 3-forms on $P$,
$\theta^i\dz\theta^j\dz\theta^k$,
$\theta^i\dz\theta^j\dz\gamma^{A/A'}$,
$\theta^i\dz\gamma^{A/A'}\dz\gamma^{B/B'}$, and
$\gamma^{A/A'}\dz\gamma^{B/B'}\dz\gamma^{C/C'}$. This brings the
relations between 
the unknown functions $u_{\mu j}, t_{\mu A}, t_{\mu A'}$ and the
curvatures $K^i_{~jkl}$. 
  
Analysing these relations step by
step we get the following:
\begin{itemize}
\item First, we consider terms at the basis forms
  $\theta^i\dz\theta^j\dz\gamma^{A/A'}$. This gives 42 conditions
  determining all the functions $u_{\mu A}$ and $u_{\mu A'}$ in terms
  of $(u_\mu)$. After solving these 42 conditions we get:
$$
\begin{aligned}
&\der u_1=(3u_4+2u_6)\gamma^1+u_5\gamma^2-2u_2\gamma^3+u_{1j}\theta^j\\
&\der u_2=-(2u_5+u_7)\gamma^1-(u_4+2u_6)\gamma^2+(u_1-u_3)\gamma^3+u_{2j}\theta^j\\
&\der u_3=u_6\gamma^1+(2u_5+3u_7)\gamma^2+2u_2\gamma^3+u_{3j}\theta^j\\
&\der u_4=-3u_1\gamma^1+3u_5\gamma^3+u_{4j}\theta^j\\
&\der u_5=2u_2\gamma^1-u_1\gamma^2+(2u_6-u_4)\gamma^3+u_{5j}\theta^j\\
&\der u_6=-u_3\gamma^1+2u_2\gamma^2+(u_7-2u_5)\gamma^3+u_{6j}\theta^j\\
&\der u_7=-3u_3\gamma^2-3u_6\gamma^3+u_{7j}\theta^j.
\end{aligned}
$$  
\item Second, the terms at the basis forms
  $\theta^i\dz\theta^j\dz\theta^k$ when equated to zero, can be split
  into two types of equations. The first type is obtained by
  eliminating the curvatures $K^i_{~jkl}$ from the full set. This
  yields a system of linear equations for the unknowns $t_{\mu j}$,
  whose only solution is $u_{\mu j}=0.$
After these conditions are imposed the second type of equations,
involves the curvatures $K^i_{~jkl}$ only in a linear fashion. It has a
unique solution for the curvatures, which explicitly is given by
(\ref{curveu})-(\ref{curveuu}).
\item Third, after imposing the conditions described above,
  all the other terms in $\der^2\theta^i$ are automatically zero. 
\end{itemize}  
This proves the proposition. 
\end{proof}
The next proposition determines the derivatives of the unknown $k_1$.
\begin{proposition}
The second Bianchi identities $\der^2\gamma^A\equiv 0\equiv \der^2\gamma^{A'}$,
$A,A'=1,2,3$, imply that
\be\der k_1=-2k_3\gamma^2+2k_2\gamma^3.\label{rio}\ee
\end{proposition}
\begin{proof}
To prove this we write $\der k_1$ in the most general form
$$
\der k_1=k_{1i}\theta^i+k_{1A}\gamma^A+k_{1A'}\gamma^{A'},$$
and consider the terms $\theta^i\dz\theta^j\dz\gamma^{A/A'}$ in the
decomposition of $\der^2\gamma^{A/A'}$. This immediately yields:
$$k_{1A'}= 0, \quad\forall A'=1,2,3,$$
and 
$$k_{11}=0,\quad k_{12}=-2k_3,\quad{\rm and}\quad k_{13}=2k_2.$$
Eliminating $u_\mu$s from the equations implied by equating to zero
the coefficients at the terms $theta^i\dz\theta^j\dz\theta^k$ in
$\der^2\gamma^{A/A'}\equiv 0$, shows that all the remaining
coefficients $k_{1i}$ in $\der k_1$ must also vanish
$$k_{1i}= 0,\quad\forall i=1,2,\dots 9.$$
This finishes the proof. 
\end{proof}

The lack of the $\theta^i$ terms on the right hand sides of equations
(\ref{toreu}) and (\ref{rio}) proves that the functions $u_\mu$ and
$k_1$, and as a consequence the functions $k_2,\dots,k_7$, are
\emph{constant} along the base manifold $M$. They depend only on the
fiber coordinates. Moreover, since only $\gamma^A$s appear on the
right hand sides of these equations, they only depend on the fiber
coordinates associated with $\sog(3)_L$. This means that there exists
a $\sog(3)_L$ gauge in which all the functions $u_\mu$,
$k_1,\dots,k_7$ are constant. This is the same to say that 
there exists a subbundle $\mathcal G$ of $P$, with fibers of at least
as large as $\sog(3)_R$, on which we have 
$$\der u_\mu=0=\der k_1=\dots=\der k_7.$$    
To see the examples of such solutions we look at the fourth and the
seventh of the equations (\ref{toreu}). Since we want $\der u_4=\der
u_7=0$,
we obtain that:
$$u_1\gamma^1=u_5\gamma^3\quad{\rm}\quad u_3\gamma^2=-u_6\gamma^3.$$
Now, assuming that $u_1\neq 0\neq u_3,$
we solve it for $\gamma^1$ and $\gamma^2$, obtaining:
\be\gamma^1=\frac{u_5}{u_1}\gamma^3,\quad{\rm and}\quad
\gamma^2=-\frac{u_6}{u_3}\gamma^3.\label{ld}\ee
Thus these equations show that we have reduced our original manifold
$P$ to its 13-dimensional submanifold $\mathcal G$ on which the forms
$\gamma^1$ and $\gamma^2$ become dependent on $\gamma^3$. On
this manifold we further want that $\der u_\mu=0$ for all
$\mu=1,2,\dots 7$. Inserting (\ref{ld}) into the right hand sides
of equations (\ref{toreu}) for $\der u_1$, $\der u_2$, $\der u_3$, $\der u_5$, $\der u_6$,
and equating the result to zero, we obtain the five equations:
$$
\begin{aligned}
&2 u_1 u_2 u_3-3 u_3 u_4 u_5+u_1 u_5 u_6-2 u_3 u_5 u_6=0\\
&u_1^2 u_3-u_1 u_3^2-2 u_3 u_5^2+u_1 u_4 u_6+2 u_1 u_6^2-u_3 u_5 u_7=0\\
&2 u_1 u_2 u_3-2 u_1 u_5 u_6+u_3 u_5 u_6-3 u_1 u_6 u_7=0\\
&u_1 u_3 u_4-2 u_2 u_3 u_5-u_1^2 u_6-2 u_1 u_3 u_6=0\\
&2 u_1 u_3 u_5+u_3^2 u_5+2 u_1 u_2 u_6-u_1 u_3 u_7=0.
\end{aligned}
$$ 
A particular solution is given by:
\be\begin{aligned}
&u_2=\frac{u_6\sqrt{4u_1+u_3}\sqrt{u_1u_3^2-u_3^3+u_1u_6^2+4u_3u_6^2}}{3u_6^2-u_3^2}\\
&u_4=\frac{u_6(u_1u_6^2-3u_1u_3^2-2u_3u_6^2)}{u_3(3u_6^2-u_3^2)}\\
&u_5=-\frac{u_1\sqrt{u_1u_3^2-u_3^3+u_1u_6^2+4u_3u_6^2}}{u_3\sqrt{4u_1+u_3}}\\
&u_7=\frac{(2u_1u_3^3+u_3^3+2u_1u_6^2-u_3u_6^2)\sqrt{u_1u_3^2-u_3^3+u_1u_6^2+4u_3u_6^2}}{u_3(3u_6^2-u_3^2)\sqrt{4u_1+u_3}}.
\end{aligned}\label{ldu}\ee
Of course we restrict the range of the free real  torsion parameters $u_1$,
$u_3$ and $u_6$, so that $u_2$, $u_4$, $u_5$ and $u_7$ are real and
finite! This happens e.g. for $-1<\tfrac{4u_3}{u_1}<4$, $u_6\neq
\pm\sqrt{\tfrac13}u_3\neq 0$.

This solution is compatible with the structure equations
$$\begin{aligned}
\der\gamma^1=-\gamma^2\dz\gamma^3+\kappa^1\\
\der\gamma^2=-\gamma^3\dz\gamma^1+\kappa^2
\end{aligned}$$
having $\kappa^1,\kappa^2$ and $\kappa^3$ as in (\ref{curveu}), and
with 
$\der k_1=0$ if and only if 
\be 
k_1=\frac{4(u_1u_3^2+u_1u_6^2+u_3u_6^2)^2(u_1u_3^2-u_3^3+u_1u_6^2+4u_3u_6^2)}{u_3^2(4u_1+u_3)(u_3^2-3u_6^2)^2}.\label{lkd}\ee
This leads to the following proposition.
\begin{proposition}
Assume that the forms $(\theta^i,\gamma^3,\gamma^{A'})$ satisfy the
equations for $\der\theta^i$, $\der\gamma^3$, and $\der\gamma^{A'}$ as
in the system
(\ref{1str})-(\ref{2str}), (\ref{mnb}) with  
\begin{itemize}
\item the forms $\gamma^1$ and
$\gamma^2$ given by (\ref{ld}),
\item the coefficients $u_1$, $u_3$ and $u_6$ being
constants, 
\item the coefficients $u_2$, $u_4$, $u_5$ and $u_7$
given by (\ref{ldu}),
\item the curvatures $\kappa^1$, $\kappa^{A'}$ given by
(\ref{curveu})-(\ref{curveuu}) and (\ref{lkd}). 
\end{itemize}
Then 
\begin{itemize}
\item the equations for $\der\gamma^1$ and $\der\gamma^2$ in the system
(\ref{1str})-(\ref{2str}), (\ref{mnb}) are automatically satisfied,
  and
\item the Bianchi identities
  $\der^2\theta^i=\der^2\gamma^3=\der^2\gamma^{A'}=0$ are also
  automatically satisfied.
\end{itemize}
In such a case the manifold on which the forms
$(\theta^i,\gamma^3,\gamma^{A'})$ are defined becomes a 
13-dimensional Lie group ${\mathcal G}^{13}$, with the forms
$(\theta^i,\gamma^3,\gamma^{A'})$ being its Maurer-Cartan forms. The Lie
group ${\mathcal G}^{13}$ is a subbundle of the 
bundle $\sog(3)\times\sog(3)\to P\to M^9$, so that the manifold $M^9$ is a
homogeneous space $M^9={\mathcal G}^{13}/H$, with $H$ being a certain 4-dimensional subgroup of ${\mathcal G}^{13}$ containing
$\sog(3)_R$. The nearly integrable $\sog(3)\times\sog(3)$ structure
$(g,\ten,\omega)$ on
$M^9$ is given by $\theta^i$s and the formulae (\ref{ty}).

For all of these geometries the metric $g$ is conformally
\emph{non}-flat and \emph{not} locally symmetric. 
The Ricci tensors of the Levi-Civita
connection $\lc$, of the characteristic connection $\Gamma$, and of the
$\soa(3)_L$ part $\gp$ of $\Gamma$ have all two 
distinct eigenvalues. 

The $\soa(3)_R$ part $\gm$ of the characteristic connection $\Gamma$  is
\emph{Einstein}. 

Explicitly, in the adapted coframe $(\theta^i)$ in which the structure
equations read as in (\ref{131}) and in which the structural tensors
$g,\ten,\omega$ are given by (\ref{ty}), we have:
\begin{itemize}
\item The eigenvalues of the Levi-Civita connection Ricci tensor read:
$$\Big(45 s,45 s,45 s,55 s,55 s,55 s,55 s,55 s,55 s\Big),$$
where 
$$s=\frac{(u_1u_3^2+u_1 u_6^2+u_3 u_6^2)^3}{u_3^2(4u_1+u_3)(u_3^2-3u_6^2)^2}.$$
The Ricci scalar is equal to $465 s.$ The Levi-Civita connection is
never Ricci flat, because the equation $u_1u_3^2+u_1 u_6^2+u_3
u_6^2=0$ contradicts the reality of $u_2$, $u_5$ and $u_7$. 
\item The $\soa(3)_L$ part $\gp$ of  $\Gamma$ has
  the curvature
$\kp=\kappa^Ae_A,$
with 
$$
\begin{aligned}
\kappa^1&=\\
&\frac{4(u_1u_3^2+u_1u_6^2+u_3u_6^2)^2(u_1u_3^2-u_3^3+u_1u_6^2+4u_3u_6^2)}{u_3^2(4u_1+u_3)(u_3^2-3u_6^2)^2}\kappa^1_0+\\
&\frac{4u_6(u_1u_3^2+u_1u_6^2+u_3u_6^2)^2\sqrt{u_1u_3^2-u_3^3+u_1u_6^2+4u_3u_6^2}}{u_3^2\sqrt{4u_1+u_3}(u_3^2-3u_6^2)^2}\kappa^2_0-
\\
&\frac{4(u_1u_3^2+u_1u_6^2+u_3u_6^2)^2\sqrt{u_1u_3^2-u_3^3+u_1u_6^2+4u_3u_6^2}}{u_3\sqrt{4u_1+u_3}(u_3^2-3u_6^2)^2}\kappa^3_0,
\end{aligned}
$$
$$
\begin{aligned}
\kappa^2&=\\
&\frac{4u_6(u_1u_3^2+u_1u_6^2+u_3u_6^2)^2\sqrt{u_1u_3^2-u_3^3+u_1u_6^2+4u_3u_6^2}}{u_3^2\sqrt{4u_1+u_3}(u_3^2-3u_6^2)^2}\kappa^1_0+\\
&\frac{4u_6^2(u_1u_3^2+u_1u_6^2+u_3u_6^2)^2}{u_3^2(u_3^2-3u_6^2)^2}\kappa^2_0-\frac{4u_6(u_1u_3^2+u_1u_6^2+u_3u_6^2)^2}{u_3(u_3^2-3u_6^2)^2}\kappa^3_0\\&
\\
\kappa^3&=\\
&-\frac{4(u_1u_3^2+u_1u_6^2+u_3u_6^2)^2\sqrt{u_1u_3^2-u_3^3+u_1u_6^2+4u_3u_6^2}}{u_3\sqrt{4u_1+u_3}(u_3^2-3u_6^2)^2}\kappa^1_0-\\
&\frac{4u_6(u_1u_3^2+u_1u_6^2+u_3u_6^2)^2}{u_3(u_3^2-3u_6^2)^2}\kappa^2_0+\frac{4(u_1u_3^2+u_1u_6^2+u_3u_6^2)^2}{(u_3^2-3u_6^2)^2}\kappa^3_0.
\end{aligned}
$$
and the matrices $e_A=(e_A\phantom{}^i_j)$ given by
(\ref{dwie}). It has the Ricci tensor $\Rp_{ij}$ with two different
eigenvalues
$$\Big(0,0,0,20 s,20 s,20 s,20 s,20 s,20 s\Big),$$
with the Ricci scalar equal to $120s$. 
\item The $\soa(3)_R$ part $\gm$ of  $\Gamma$  has
  the curvature
$\km=15s\kappa_0^{A'}e_{A'},$
where as before the matrices $e_{A'}=(e_{A'}\phantom{}^i_j)$ are given by
(\ref{dwie}). Its Ricci tensor is \emph{Einstein},
$\Rm_{ij}=30sg_{ij},$
and has Ricci scalar equal to $270s$. 
\item The characteristic connection $\Gamma=\gp+\gm$ has curvature
$$\Omega=\kp+\km=\kappa^Ae_A+15s\kappa_0^{A'}e_{A'}$$
and the Ricci
  tensor with eigenvalues:
$$\Big(30 s,30 s,30 s,50 s,50 s,50 s,50 s,50 s,50 s\Big).$$
 \end{itemize}
\end{proposition}
The examples of nearly integrable
$\sog(3)\times\sog(3)$ geometries with torsion of the
characteristic connection in $V_{[0,6]}$ described by this proposition
have quite similar features to the nearly integrable
$\sog(3)\times\sog(3)$ geometries with torsion in $V_{[0,2]}$. In
particular, if any of these geometries has curvature $\kp\equiv 0$
then it must be flat, and torsion free. 

It turns out however that there is another branch of nearly integrable
$\sog(3)\times\sog(3)$ geometries with torsion of their characteristic
connections in $V_{[0,6]}$ for which $\kp\equiv 0$ does not imply
neither vanishing torsion nor vanishing of $\km$. Below we present
these examples.   

Assuming that $$\kp\equiv 0$$ is the same as to assume that 
$k_1=k_2=k_3=k_4=k_5=k_6=0$. (Compare with the first three equations (\ref{curveu})). But since
$\kp\equiv 0$ is the condition for the connection $\gp$ to be flat, in
such a situation we can use a gauge in which $\gp\equiv 0$. This
condition means that the system (\ref{1str})-(\ref{2str}), (\ref{mnb})
reduces form $P$ to a 12-dimensional ${\mathcal G}^{12}$ manifold on which
$$\theta^{10}\equiv\theta^{11}\equiv\theta^{12}\equiv 0.$$ 
Having these conditions and the requirement that $T\in V_{[0,6]}$ implies, via
(\ref{toreu}), that all $u_\mu$ are \emph{constants}. The rest of the
equations $\der^2\theta^i\equiv 0$ imply finally that:
$$
\begin{aligned}
&2 u_2 u_4+2 u_1 u_5-u_3 u_5+u_1 u_7=0\\
&u_3 u_4-u_1 u_6+2 u_3 u_6+2 u_2 u_7=0\\
&2 u_1 u_2+2 u_2 u_3-2 u_4 u_5-3 u_5 u_6-u_4 u_7-2 u_6
  u_7=0\\
&2 u_1^2+2 u_1 u_3+2 u_4^2+4 u_4 u_6+2 u_5 u_7=0\\
&2 u_1 u_3+2 u_3^2+2 u_4 u_6+4 u_5 u_7+2 u_7^2=0\\
&2 u_1^2-2 u_3^2+2 u_4^2+2 u_4 u_6-2 u_5 u_7-2 u_7^2=0.
\end{aligned} 
$$
We have found 6 different particular solutions to these equations. These are:
\begin{enumerate}
\item 
$\begin{aligned}
&u_2=\frac{(u_1-2u_3)u_7^2-(2u_1-u_3)\Big((u_1-2u_3)(u_1+u_3)+u_4^2\Big)}{6u_4u_7}\\[3 pt]
&u_5=\frac{(u_1-2u_3)(u_1+u_3)+u_4^2-2u_7^2}{3 u_7}, \quad u_6=\frac{-(2u_1-u_3)(u_1+u_3)-2u_4^2+u_7^2}{3u_4};
\end{aligned}$\\[3 pt]
\item
$\begin{aligned}
&u_2=\mp\frac{u_5\sqrt{9u_3^2-4u_4^2}}{2u_4},\quad
  u_6=\mp\frac{u_3(\pm3u_3+\sqrt{9u_3^2-4u_4^2})}{2u_4}\\[3 pt] 
  &
  u_1=\tfrac12(u_3\pm\sqrt{9u_3^2-4u_4^2}),\quad u_7=0;\end{aligned}$
\\[3 pt]
\item
$\begin{aligned}
&u_2=\mp\frac{u_6(\pm 9u_3^2+\sqrt{9u_3^2+8u_7^2)}}{8u_7},\quad
  u_5=\frac{-4u_7^2\pm u_3(\mp3u_3+\sqrt{9u_3^2+8u_7^2})}{8u_7}\\&
  u_1=\tfrac14(-u_3\mp\sqrt{9u_3^2+8u_7^2}),\quad u_4=0;\end{aligned}$
\\[3 pt]
\item
$\begin{aligned}
&u_1=u_3=u_4=u_5=u_7=0;\end{aligned}$
\\[3 pt]
\item
$\begin{aligned}
&u_1=-u_3,\quad  u_4=u_5=u_6=u_7=0;\end{aligned}$
\\[3 pt]
\item
$\begin{aligned}
&u_1=u_3=u_4=u_6=u_7=0.\end{aligned}$
\end{enumerate}
It follows that for all of these 6 solutions we have
$\der^2\theta^i\equiv 0$ and $\der\gamma^{A'}\equiv 0$, automatically
for all $i=1,2,\dots 9$ and for all $A'=1,2,3$. Thus each of these 6
solutions defines a nearly integrable $\sog(3)\times\sog(3)$ geometry
$(M^9,g,\ten,\omega)$ having the torsion of the characteristic
connection in $V_{[0,6]}$ and the vanishing curvature $\kp$ of 
 $\gp$. It turns out that all the six solutions have the
same qualitative behaviour of the curvatures of $\lc$, $\Gamma$, $\gp$
and $\gm$. The properties of
the curvatures of the geometries corresponding to these six solutions
are summarized in the theorem below.    
\begin{theorem}\label{v06}
All nearly integrable $\sog(3)\times\sog(3)$ geometries
$(M^9,g,\ten,\omega)$  corresponding to any solution (1)-(6) above have 
\begin{itemize}
\item torsion of the characteristic connection $\Gamma$ in
$V_{[0,6]}\subset\bgw^3\bbR^9$
\item vanishing curvature $\kp$ of the
$\soa(3)_L$ part of  $\Gamma$,  i.e. $\kp\equiv 0$
\item the curvature $\Omega$ of the characteristic connection $\Gamma$
  equal to 
$$\Omega~\equiv ~\km=\tfrac{1}{36} \| T \|^2\kappa^{A'}_0e_{A'},$$
 where $\|T \|^2$ is the square norm of the torsion $T$ of  $\Gamma$:
$$\|T \|^2=T_{ijk}T^{ijk}=36k_7=36(2u_1^2+u_2^2+u_1u_3+2u_4^2+u_5^2+3u_4u_6+u_6^2+u_5u_7)$$
with $u_\mu$ being constants and satisfying one of 
(1)-(6).
\end{itemize}
All these geometries $(M^9,g,\ten,\omega)$ are locally homogeneous spaces
$M^9={\mathcal G}^{12}/H$, where ${\mathcal G}^{12}$ is a 12-dimensional
  symmetry group of $(M^9,g,\ten,\omega)$ and $H$ is its 3-dimensional subgroup
isomorphic to $\sog(3)$, $H=\sog(3)_R$. The metric $g$, the tensor
$\ten$ and the form $\omega$ defining a nearly integrable
$\sog(3)\times\sog(3)$ geometry on $M^9$ are given by formulae
(\ref{ty}), in terms of the forms
$(\theta^i,\gamma^{A}\equiv 0,\gamma^{A'})$ satisfying
(\ref{1str})-(\ref{2str}), (\ref{mnb}),
(\ref{curveu})-(\ref{curveuu}), and one of (1)-(6), with $u_\mu$ being
constants.
\begin{itemize}
\item In the basis $(\theta^i,\gamma^{A'})$ the Killing form for the group ${\mathcal G}^{12}$ reads:
$$Kil~=~-8~{\rm
  diag}\Big(k_7,k_7,k_7,k_7,k_7,k_7,k_7,k_7,k_7,1,1,1\Big).$$
\item If $k_7\neq 0$ the Riemannian manifold
  $(M^9={\mathcal G}^{12}/\sog(3)_R,  g)$ is not locally symmetric. If $k_7=0$ the
  solutions have flat characteristic connection, $\Omega\equiv 0$, and
  in such a case $(M^9={\mathcal G}^{12}/\sog(3)_R,  g)$ is a locally
  symmetric Riemannian manifold. 
\item For every value of $k_7$ the metric
  is Einstein, $\rilc = 3k_7g.$ It is not conformally flat unless the torsion is
  zero, $(u_1,u_2,\dots,u_7)=0$. 
\item Also the $\sog(3)_R$ part $\gm$ of the characteristic
  connection is always Einstein, $\Rm_{ij}=2k_7g_{ij}.$
It is flat, $\km\equiv 0$, if and only if $k_7= 0$. 
\end{itemize}   
\end{theorem}
It is a remarkable
fact that both the Levi-Civita connection $\lc$ and the characteristic
connection $\Gamma$ are 
Einstein and (generically) Ricci non flat for all the geometries
$(M^9,g,\ten,\omega)$ described by the theorem. 
Moreover although the metric $g$ is  not conformally flat, the $\sog(3)_L$
part $\gp$ of $\Gamma$ is flat. This makes these
geometries similar to the selfdual Riemannian geometries in dimension
four.
\subsection{Analogs of selfduality; examples with torsion in
  $V_{[0,2]}\oplus V_{[0,6]}$ } 
The examples described by the Theorem \ref{v06} raise the question if
there are other nearly integrable $\sog(3)\times\sog(3)$ geometries
$(M^9,g,\ten,\omega)$ in dimension nine for which the $\soa(3)_L$ part
$\gp$ of the characteristic connection $\Gamma$ is flat, $\kp\equiv
0,$ and for which the $\soa(3)_R$ part $\gm$, is not flat, $\km\neq
0.$ 

 In the following the nearly integrable $\sog(3)\times\sog(3)$ geometries
$(M^9,g,\ten,\omega)$  with these two properties,  $\kp\equiv
0$ and $\km\neq
0,$ will be called \emph{analogs of
selduality}.

The problem of finding all such structures is a difficult one. To
generalize solutions of Theorem \ref{v06}, on top of the analogs of
selfduality conditions, we will assume in addition
that the torsion $T$ of the characteristic connection $\Gamma$ is
restricted from $\bgw^3\bbR^9$ to 
$V_{[0,2]}\oplus V_{[0,6]}$.
In this section we will find all such structures. 

We first have an analog of Propositions
\ref{prtor7} and Remark \ref{prtor1}:
\begin{proposition}\label{prtor10}
In an adapted coframe $(\theta^i)$ the $V_{[0,2]}\oplus V_{[0,6]}$ torsion 
of a characteristic connection of a nearly integrable 
geometry $(M^9,g,\ten,\omega)$ reads:
\be
\begin{aligned}
T^1&=-t_1\lambda_0^{9'}+t_2\lambda_0^{6'}+\tfrac13 t_3(5\kappa_0^{3'}-4\lambda_0^{3'}+2\lambda_0^{12'})+\\&u_1(-\lambda_0^{3'}+\lambda_0^{12'})-u_2\lambda_0^{15'}-u_3\lambda_0^{3'}-u_4\lambda_0^{6'}-u_5\lambda_0^{9'}-u_6\lambda_0^{6'}-u_7\lambda_0^{9'}\\
&\\
T^2&=t_1\lambda_0^{8'}-t_2\lambda_0^{5'}+\tfrac13 t_3(-5\kappa_0^{2'}+4\lambda_0^{2'}-2\lambda_0^{11'})+\\&u_1(\lambda_0^{2'}-\lambda_0^{11'})+u_2\lambda_0^{14'}+u_3\lambda_0^{2'}+u_4\lambda_0^{5'}+u_5\lambda_0^{8'}+u_6\lambda_0^{5'}+u_7\lambda_0^{8'}\\
&\\
T^3&=-t_1\lambda_0^{7'}+t_2\lambda_0^{4'}+\tfrac13 t_3(5\kappa_0^{1'}-4\lambda_0^{1'}+2\lambda_0^{10'})+\\&u_1(-\lambda_0^{1'}+\lambda_0^{10'})-u_2\lambda_0^{13'}-u_3\lambda_0^{1'}-u_4\lambda_0^{4'}-u_5\lambda_0^{7'}-u_6\lambda_0^{4'}-u_7\lambda_0^{7'}\\
&\\
T^4&=-t_1\lambda_0^{15'}+\tfrac13
  t_2(-5\kappa_0^{3'}-2\lambda_0^{3'}+4\lambda_0^{12'})-t_3\lambda_0^{6'}+\\&u_1\lambda_0^{6'}-u_2\lambda_0^{9'}+u_4(-\lambda_0^{3'}+\lambda_0^{12'})+u_5\lambda_0^{15'}-u_6\lambda_0^{3'}\\
&\\
T^5&=t_1\lambda_0^{14'}+\tfrac13
  t_2(5\kappa_0^{2'}+2\lambda_0^{2'}-4\lambda_0^{11'})+t_3\lambda_0^{5'}-\\&u_1\lambda_0^{5'}+u_2\lambda_0^{8'}+u_4(\lambda_0^{2'}-\lambda_0^{11'})-u_5\lambda_0^{14'}+u_6\lambda_0^{2'}\\
&\\
T^6&=-t_1\lambda_0^{13'}+\tfrac13 t_2(-5\kappa_0^{1'}-2\lambda_0^{1'}+4\lambda_0^{10'})-t_3\lambda_0^{4'}+\\&u_1\lambda_0^{4'}-u_2\lambda_0^{7'}+u_4(-\lambda_0^{1'}+\lambda_0^{10'})+u_5\lambda_0^{13'}-u_6\lambda_0^{1'}\\
&\\
T^7&=\tfrac13
  t_1(5\kappa_0^{3'}+2\lambda_0^{3'}+2\lambda_0^{12'})+t_2\lambda_0^{15'}-t_3\lambda_0^{9'}-\\&u_2\lambda_0^{6'}+u_3\lambda_0^{9'}+u_5(-\lambda_0^{3'}+\lambda_0^{12'})+u_6\lambda_0^{15'}-u_7\lambda_0^{3'}\\
&\\
T^8&=-\tfrac13
  t_1(5\kappa_0^{2'}+2\lambda_0^{2'}+2\lambda_0^{11'})-t_2\lambda_0^{14'}+t_3\lambda_0^{8'}+\\&u_2\lambda_0^{5'}-u_3\lambda_0^{8'}+u_5(\lambda_0^{2'}-\lambda_0^{11'})-u_6\lambda_0^{14'}+u_7\lambda_0^{2'}\\
&\\
T^9&=\tfrac13 t_1(5\kappa_0^{1'}+2\lambda_0^{1'}+2\lambda_0^{10'})+t_2\lambda_0^{13'}-t_3\lambda_0^{7'}-\\&u_2\lambda_0^{4'}+u_3\lambda_0^{7'}+u_5(-\lambda_0^{1'}+\lambda_0^{10'})+u_6\lambda_0^{13'}-u_7\lambda_0^{1'},
\end{aligned}\label{mnb1}
\ee 
where $(t_1,t_2,t_3,u_1,u_2,u_3,u_4,u_5,u_6,u_7)$ are the ten independent components of the
torsion $T$, and $(\lambda_0^{\mu'})$, $\mu'=1,2,\dots, 15$, is a basis of 2-forms in
$V_{[4,2]}$ as defined in (\ref{lambda0p}). Note that if all $u_\mu$s
are equal zero the torsion $T\in V_{[0,2]}$, and if all $t_A$s are
equal zero $T\in V_{[0,6]}$.  
\end{proposition}      
We want to construct nearly integrable $\sog(3)\times\sog(3)$
structures with torsion in $V_{[0,2]}\oplus V_{[0,6]}$, and with 
$\kp\equiv 0$. All of them, in an adapted
coframe, are therefore described by the system (\ref{1str})-(\ref{2str}), (\ref{mnb1}), with
$\kappa^A\equiv 0$. This enables us to reduce the system 
from $P\to M^9$ to a 12 dimensional subbundle of $P$ on which 
$$\theta^{10}\equiv\theta^{11}\equiv\theta^{12}\equiv 0. $$ 

The procedure of analysing such a reduced system is completely the
same as the procedure leading to solutions described by the Theorem
\ref{v06}. We therefore only state the result.

\begin{theorem}
All nearly integrable $\sog(3)\times\sog(3)$ geometries
$(M^9,g,\ten,\omega)$, which have torsion $T$ of the
characteristic connection $\Gamma$ in $V_{[0,2]}\oplus V_{[0,6]}$, and the
curvature $\kp$ of the $\soa(3)_L$-part $\gp$ of  $\Gamma$ vanishing, $\kp\equiv 0$, correspond to the system   
(\ref{1str})-(\ref{2str}), (\ref{mnb1}), with
$$\theta^{10}\equiv\theta^{11}\equiv\theta^{12}\equiv 0,\quad\quad\kappa^A\equiv 0,$$ and \emph{constant} torsion coefficients
$(t_1,t_2,t_3,u_1,u_2,u_3,u_4,u_5,u_6,u_7)$ satisfying the following
algebraic equations:
\be
\begin{aligned}
2 u_2 u_4+2 u_1 u_5&-u_3 u_5+u_1 u_7+\\&t_2 u_2+t_1 u_3-t_3 u_5-t_3
  u_7-t_1 t_3=0\\[5 pt]
u_3 u_4-u_1 u_6&+2 u_3 u_6+2
  u_2 u_7-\\&t_2 u_1-t_1 u_2-t_3 u_4-t_3 u_6+t_2 t_3=0\\[5pt]
2 u_1 u_2+2 u_2 u_3&-2 u_4 u_5-3 u_5
  u_6-u_4 u_7-2 u_6 u_7+\\&t_3 u_2+t_2 u_5-t_1 u_6-t_1 t_2=0\\[5pt]
2 u_1^2+2 u_1 u_3&+2 u_4^2+4 u_4 u_6+2 u_5 u_7-\\&2 t_1 u_7-t_3 u_1-2 t_3 u_3+t_2 u_4-t_1 u_5+2 t_2 u_6+t_1^2-t_3^2=0\\[5pt]
2 u_1 u_3+2 u_3^2&+2 u_4 u_6+4 u_5 u_7+2 u_7^2-\\&2 t_3 u_1-t_3 u_3+2 t_2 u_4-2 t_1 u_5+t_2 u_6-t_1 u_7+t_2^2-t_3^2=0\\[5pt]
2 u_1^2-2 u_3^2&+2 u_4^2+2 u_4 u_6-2 u_5 u_7-2 u_7^2+\\&t_3 u_1-t_3 u_3-t_2 u_4+t_1
  u_5+t_2 u_6-t_1 u_7+t_1^2-t_2^2=0.
\end{aligned}
\label{feq}\ee
If these equations are satisfied the metric $g$, the tensor $\ten$ and
the 4-form $\omega$ are obtained in terms of the forms $(\theta^i)$
via formulae (\ref{ty}). They descend from the 12-dimensional
subbundle $P^{12}\to M^9$ of the fiber bundle $\sog(3)\times\sog(3)\to
P\to M^9$ to $M^9$ due to the structure equations (\ref{1str}). 

If the equations (\ref{feq}) for the constants
$(t_1,t_2,t_3,u_1,u_2,u_3,u_4,u_5,u_6,u_7)$ are satisfied, then all
the integrability condions $\der^2\theta^i\equiv 0$ and
$\der\gamma^{A'}\equiv 0$, for all $\theta^i$s and $\gamma^{A'}$s appearing in the system
(\ref{1str})-(\ref{2str}), (\ref{mnb1}) are automatically
satisfied. 

The manifold $P^{12}$ is locally a 12-dimensional symmetry group
$P^{12}={\mathcal G}^{12}$ of
the so obtained $(M^9,g,\ten,\omega)$, and $M^9$ is a homogeneous space
$M^9={\mathcal G}^{12}/H$, where $H=\sog(3)_R$ is a subgroup of
${\mathcal G}^{12}$. 

The curvatures $\kappa^{A'}$ are given by
$$\kappa^{A'}=\Big(~\tfrac{1}{36} \|T \|^2-\tfrac{25}{6}(t_1^2+t_2^2+t_3^2)~\Big)~\kappa^{A'}_0, \quad\quad A'=1,2,3,$$
where 
$$\begin{array}{lcl}
\|T \|^2 &= &6(4 u_1^2+6 u_2^2+2 u_1 u_3+4 u_3^2+4 u_4^2+6
u_5^2+6 u_4 u_6+6 u_6^2+6 u_5 u_7+4 u_7^2)\\[3pt]
&&+90(t_1^2+t_2^2+t_3^2),\end{array}$$
with constants $(t_1,t_2,t_3,u_1,u_2,u_3,u_4,u_5,u_6,u_7)$  fulfilling
equations (\ref{feq}).

The torsion $T$ of  the characteristic connection  generically seats in 
$V_{[02,]}\oplus V_{[0,6]}$. It is in $V_{[0,2]}$ iff
$(u_1,u_2,u_3,u_4,u_5,u_6,u_7)=0$, and in $V_{[0,6]}$ iff
$(t_1,t_2,t_3)=0$. The square of the torsion is $\|T\|^2$ as above. 

The curvature $\Omega$ of  $\Gamma$ has vanishing
$\soa(3)_L$ part, $\kp\equiv 0$, and is equal to:
$$\Omega\equiv \km=
\Big(~\tfrac{1}{36} \| T \|^2-\tfrac{25}{6}(t_1^2+t_2^2+t_3^2)~\Big)~\kappa^{A'}_0e_{A'}.$$
The Ricci tensor of the curvature $\Omega$ of the characteristic
connection, and what is the same, the Ricci tensor of the curvature
$\km$ of
its $\soa(3)_R$-part is Einstein,
$$\Rm_{ij}=2\Big(~\tfrac{1}{36} \|T \|^2-\tfrac{25}{6}(t_1^2+t_2^2+t_3^2)~\Big)g_{ij}.$$
The metric $g$ is Einstein if and only if $t_1=t_2=t_3=0$. In such a
case the nearly integrable structures coincide with those described in
Theorem \ref{v06}. 

Generically the solutions described by this theorem have $\km\neq
0$, and as such constitute analogs of selfduality. 
\end{theorem}  
\begin{remark}
Note that although $(t_1,t_2,t_3)=0$ gives all the solutions described
in Theorem \ref{v06}, the assumption $(u_1,u_2,u_3,u_4,u_5,u_6,u_7)=0$
does not recover all the solutions with $T\in V_{[0,2]}$. The reason for
this is that here we additionally assumed that $\kp\equiv 0$, and such
solutions are possible for $T\in V_{[0,2]}$ only if $T=0$. Nonetheless
the solutions in this section are nontrivial generalizations to $T\in
V_{[0,2]}\oplus V_{[0,6]}$ of
solutions from Theorem \ref{v02} and \ref{v06}. 
\end{remark}
\begin{remark}
We emphasize that the system of equations (\ref{feq}) for the
constants $(t_1,t_2,t_3,u_1,u_2,u_3,u_4,u_5,u_6,u_7)$ can be solved
explicitly to the very end. For example, an application of a 
Mathematica command {\tt Solve[]} to the system (\ref{feq}), immediately gives 13 different solutions of these
equations. The obtained formulae are not particularly illuminating. For example
a generalization to the case of $T\in V_{[0,2]}\oplus V_{[0,6]}$ of
the solution (1) from Section \ref{sv06} is given by:
$$\begin{aligned}
u_2&=\Big(\big(2 (t_3+u_1-2 u_3) u_7^2+\\&(t_3-2 u_1+u_3) (-2
    t_2^2+(t_3+u_1-2 u_3) (t_3+2 (u_1+u_3))-3 t_2 u_4+2 u_4^2)+\\&3 t_1
    (t_3+u_1-2 u_3) u_7-2 t_1^2 (t_3+u_1-2 u_3)\big)\Big)\times\Big(
    3(t_2+2u_4)(2u_7-t_1)\Big)^{-1}\\&&\\
u_5&=\frac{(t_3+u_1-2 u_3) (t_3+2 (u_1+u_3))-(2 t_2-u_4) (t_2+2 u_4)-4 u_7^2+t_1^2}{3(2 u_7-t_1)},\\
&\\
&u_6=\frac{2 u_3 (u_3-u_1)-4 (u_1^2+u_4^2)-(t_1-2 u_7) (2 t_1+u_7)+t_2^2+t_3^2+3 t_3 u_3}{3(2u_4+t_2)}.
\end{aligned}
$$
It is a matter of cheking that this becomes a solution (1) from
Section \ref{sv06} in the limit $t_1\to 0$, $t_2\to 0$, $t_3\to 0$.

A solution of (\ref{feq}) which has \emph{no} limit when $t_1\to 0$, $t_2\to
0$, $t_3\to 0$ is given below:
$$\begin{aligned}
u_2=\frac{3 t_1 t_2-8 t_2 u_5+8 t_1 u_6+12 u_5 u_6}{20t_3},\quad
u_1=u_3=t_3,\quad u_4=-\tfrac12 t_2,\quad u_7=\tfrac12 t_1.
\end{aligned}
$$
\end{remark}
\begin{remark}
It is remarkable that we have obtained analogs of selfduality
\emph{with high number of symmetries}. We did not \emph{assume} any
symmetry conditions. The homogeneity of the structures obtained were
implied by the merely requirements that $\kp\equiv 0$ and $T\in
V_{[0,2]}\oplus V_{[0,6]}$. It would be very interesting to find
analogs of selfduality which are not locally homogeneous. If such
solutions may exist is an open question. 
\end{remark}

\bigbreak


\begin{thebibliography}{99}



\bibitem{ABF} I. Agricola, J. Becker-Bender, T. Friedrich, On the topology and the geometry of SO(3)-manifolds,  \emph{Ann. Global Anal. Geom.}  {\bf 40} (2011), 67--84. 


\bibitem{AHS} M.F. Atiyah, N.J. Hitchin and I.M. Singer, Self-duality in four dimensional Riemannian geometry, \emph{Proc. Roy. Soc. London}  {\bf A362} (1978), 425--461.

\bibitem{bobi}  M.  Bobie\'nski,  P. Nurowski, Irreducible $\sog(3)$ 
  geometries in dimension five, {\it J. Reine. angew. Math.} {\bf 605} (2007), 
  51--93.
  \bibitem{br} R. Bryant, private communication, (2011).
\bibitem{cf}  S. Chiossi,  A. Fino,  Nearly integrable
  SO(3)-structures on 5-dimensional Lie groups,  \emph{J. of Lie Theory}
  {\bf 17}  (2007), 539--562.
\bibitem{dub} B. Doubrov, private communication, (2011).
  \bibitem{GMM} H. Gluck, D. Mackenzie, F. Morgan, Volume-Minimizing Cycles in Grassmann Manifolds,  \emph{Duke Math J.} {\bf 79}  (1995), 335--404.
\bibitem{god} 
M. Godli\'nski,  P. Nurowski, $\glg(2,\bbR)$ geometries of
  ODEs,  \emph{J. Geom. Phys.}  {\bf 60} (2010), 991--1027.
\bibitem{gut} J. Gutt, Special Riemannian geometries and the magic
  square of Lie algebras, http://arxiv.org/abs/0810.2138, (2008).
\bibitem{nh} C. D. Hill, P. Nurowski, Differential equations and
  para-CR structures, \emph{Bollettino dell'Unione Matematica
    Italiana}, {\bf 9 III} (2010) 25-91. 
    \bibitem{mic} C. T  Michael, The Uniqueness of Calibrated Cycles in Grassmann Manifolds, Duke PhD thesis (1996).
\bibitem{nur0} 
P. Nurowski, Distinguished dimensions for special Riemannian geometries
\emph{J. Geom. Phys.} {\bf 58} (2008),  1148--1170.
\bibitem{nurm1} 
P. Nurowski Special Riemannian geometries modeled on distinguished symmetric
 spaces {\it  RIMS Kokyuroku} {\bf 1502}  (2006) ,  107--116, math.DG/0603663.
 \bibitem{olver} P. J. Olver, {\it Classical invariant theory},  LMS Student Texts {\bf 44}, Cambridge University Press, (1999).
\bibitem{peano} G. Peano, Formazioni invariantive delle
  corrispondenze,   \emph{Giornale di
Matematiche ad uso degli studenti delle Universit\`a  italiane} {\bf 20} (1882), 79--100.
  
  \bibitem{ST}    I. M.   Singer,  J. A. Thorpe, The curvature of 4-dimensional Einstein spaces, \emph{Global Analysis} (Papers in Honor of K. Kodaira), Univ. Tokyo Press, (1969), 355--365.


  
  
  
 \end{thebibliography}
\end{document}